\documentclass[11pt]{article}
\usepackage[utf8]{inputenc}
\usepackage{amsmath,amssymb,amscd,epsfig,bbm}
\usepackage{stmaryrd,mathabx}
\usepackage{comment}
\usepackage{color}
\usepackage[T1]{fontenc}

\usepackage{interval}

\usepackage[textsize=scriptsize,textwidth=2cm]{todonotes}
\usepackage{enumitem}
\usepackage{varwidth}
\setlist{nolistsep}
\usepackage[colorlinks]{hyperref}

\usepackage{mathtools}
\makeatletter
\DeclareRobustCommand\widecheck[1]{{\mathpalette\@widecheck{#1}}}
\def\@widecheck#1#2{%
    \setbox\z@\hbox{\m@th$#1#2$}%
    \setbox\tw@\hbox{\m@th$#1%
       \widehat{%
          \vrule\@width\z@\@height\ht\z@
          \vrule\@height\z@\@width\wd\z@}$}%
    \dp\tw@-\ht\z@
    \@tempdima\ht\z@ \advance\@tempdima2\ht\tw@ \divide\@tempdima\thr@@
    \setbox\tw@\hbox{%
       \raise\@tempdima\hbox{\scalebox{1}[-1]{\lower\@tempdima\box
\tw@}}}%
    {\ooalign{\box\tw@ \cr \box\z@}}}
\makeatother

\usepackage{amsthm}

\newtheorem{defi}{Definition}[section]
\newtheorem{prop}[defi]{Proposition}
\newtheorem{theo}[defi]{Theorem}
\newtheorem{theofr}[defi]{Théorème}
\newtheorem{conj}[defi]{Conjecture}
\newtheorem{lemm}[defi]{Lemma}
\newtheorem{lemmfr}[defi]{Lemme}
\newtheorem{coro}[defi]{Corollary}

\theoremstyle{definition}
\newtheorem{rema}[defi]{Remark}
\newtheorem{exem}[defi]{Example}
\newtheorem{exems}[defi]{Examples}

\newcommand{\bdefi}{\begin{defi}}
\newcommand{\edefi}{\end{defi}}
\newcommand{\bprop}{\begin{prop}}
\newcommand{\eprop}{\end{prop}}
\newcommand{\btheo}{\begin{theo}}
\newcommand{\etheo}{\end{theo}}
\newcommand{\btheofr}{\begin{theofr}}
\newcommand{\etheofr}{\end{theofr}}
\newcommand{\blemm}{\begin{lemm}}
\newcommand{\elemm}{\end{lemm}}
\newcommand{\blemmfr}{\begin{lemmfr}}
\newcommand{\elemmfr}{\end{lemmfr}}
\newcommand{\brema}{\begin{rema}}
\newcommand{\erema}{\end{rema}}
\newcommand{\bexer}{\begin{exem}}
\newcommand{\eexer}{\end{exem}}
\newcommand{\bexems}{\begin{exems}}
\newcommand{\eexems}{\end{exems}}
\newcommand{\bconj}{\begin{conj}}
\newcommand{\econj}{\end{conj}}
\newcommand{\bcoro}{\begin{coro}}
\newcommand{\ecoro}{\end{coro}}
\newcommand{\dem}{\noindent{\bf Proof. }}

\usepackage{mathrsfs}
\renewcommand\mathcal{\mathscr}

\newcommand{\A}{{\cal A}}

\newcommand{\D}{{\cal D}}

\newcommand{\G}{{\cal G}}

\newcommand{\N}{{\cal N}}

\newcommand{\T}{{\cal T}}

\newcommand{\V}{{\cal V}}

\newcommand{\maths}[1]{{\mathbb #1}}  

\newcommand{\CC}{\maths{C}}
\newcommand{\DD}{\maths{D}}
\newcommand{\FF}{\maths{F}}
\newcommand{\HH}{\maths{H}}

\newcommand{\NN}{\maths{N}}

\newcommand{\PP}{\maths{P}}

\newcommand{\RR}{\maths{R}}

\newcommand{\XX}{\maths{X}}

\newcommand{\ZZ}{\maths{Z}}

\newcommand{\weakstar}{\overset{*}\rightharpoonup}
\newcommand{\ra}{\rightarrow}
\newcommand{\bs}{\backslash}

\newcommand{\wt}[1]{{\widetilde{#1}}}
\newcommand{\wh}[1]{{\widehat{#1}}}

\newcommand{\ga}{\gamma}
\newcommand{\Ga}{\Gamma}
\newcommand{\ssm}{\!\smallsetminus\!}

\newcommand{\cqfd}{\hfill$\Box$}

\newcommand{\bigO}{\operatorname{O}}

\newcommand{\card}{{\operatorname{Card}}}
\newcommand{\CAT}{\operatorname{CAT}}

\newcommand{\covol}{\operatorname{covol}}

\newcommand{\diam}{{\operatorname{diam}}}
\newcommand{\diag}{{\operatorname{diag}}}

\newcommand{\Div}{\operatorname{Div}}

\newcommand{\End}{\operatorname{End}}

\newcommand{\gengeod}{\operatorname{\widecheck{\G\,}\!\!}}

\renewcommand{\Im}{{\operatorname{Im}}}
\newcommand{\Isom}{\operatorname{Isom}}
\newcommand{\Leb}{\operatorname{Leb}}

\newcommand{\mBM}{m_{\rm BM}}

\newcommand{\Par}{\operatorname{Par}}

\newcommand\Perp{\operatorname{Perp}}

\renewcommand{\Re}{{\operatorname{Re}}}

\newcommand{\wtmBM}{\wt m_{\rm BM}}

\newcommand{\hdr}{{\HH}^2_\RR}

\newcommand{\PSL}{\operatorname{PSL}}

\newcommand{\PSLR}{\operatorname{PSL}_{2}(\RR)}

\newcommand{\PGL}{\operatorname{PGL}}

\newcommand{\flow}[1]{{{\tt g}^{#1}}}  

\newcommand\normalout{\partial^1_{+}}
\newcommand\normalin{\partial^1_{-}}
\newcommand\normalpm{\partial^1_{\pm}}

\newcounter{fig}

\usepackage{ulem}

\title{Equidistribution of divergent geodesics \\ in negative
  curvature} \author{Jouni Parkkonen, Frédéric Paulin, Rafael
  Sayous} \date{\today}

\begin{document}
\bibliographystyle{alphas}
\maketitle

\begin{abstract} 
In the unit tangent bundle of noncompact finite volume negatively
curved Riemannian manifolds, we prove the equidistribution towards the
measure of maximal entropy for the geodesic flow of the Lebesgue
measure along the divergent geodesic flow orbits, as their complexity
tends to infinity. We prove the analogous result for geometrically
finite tree quotients, where the equidistribution takes place in the
quotient space of geodesic lines towards the Bowen-Margulis measure.
\footnote{{\bf Keywords:} equidistribution, counting, divergent
geodesic, geodesic flow, negative curvature, trees, equilibrium
states.~~ {\bf AMS codes:} 37D40, 53C22, 20E08, 37D35, 37A25.}
\end{abstract}

\section{Introduction}
\label{sec:intro}

Let $M$ be a finite volume complete connected Riemannian good orbifold 
with dimension at least $2$ and pinched negative curvature at most
$-1$. A (locally) geodesic line $\ell$ in $M$ (or its unit tangent
vector $\dot \ell(0)\,$) is {\it divergent} if the map $\ell:\RR\ra M$
is a proper map. In this paper, we define a natural type and a natural
complexity of the divergent geodesic lines, and study their counting
and equidistribution properties in the unit tangent bundle $T^1M$ of
$M$ as their complexity tends to $+\infty$ in any given type.

The study of counting and equidistribution properties of similarly
defined divergent flats in finite volume arithmetic nonpositively
curved locally symmetric spaces has produced many works, see for
instance \cite{TomWei03, Weiss06, DavSha18, ShaZhe19, DavSha20,
  SolTam23, DanPauSay25}. See also \cite{Pollicott21} for a specific
counting problem of divergent geodesics in geometrically finite
Kleinian manifolds. Curiously, the problem of the equidistribution of
divergent geodesics in general negatively curved manifolds (in
particular with variable curvature), does not seem to have been
studied so far. One possibility to go around the noncompactness of
$M$, as in \cite{DavSha18, ShaZhe19, DavSha20}, could be to work in
the projective space of locally finite positive Borel measures on
$T^1M$ for the quotient topology of the weak-star topology. But as in
\cite{DanPauSay25}, we prefer a more precise result, exhibiting the
precise scaling factor, that immediately implies the projective
convergence.

Let us introduce some definitions and notations before stating our
main result. We denote by $\pi:T^1M\ra M$ the footpoint projection, by
$(\flow{t})_{t\in\RR}$ the geodesic flow on $T^1M$ and by $h_M$ the
topological entropy of $(\flow{t})_{t\in\RR}$.  We denote by $\End(M)$
the finite set of ends\footnote{See for instance
\cite[Sect.\ I.8.27]{BriHae99} for the definition.} of the locally
compact space $M$. We fix a family $(\V_e)_{e\in \End(M)}$ of closed
Margulis neighborhoods of the ends of $M$ with pairwise disjoint
interiors, see Section \ref{sec:divgeod} for a definition. For every
$A\geq 0$, we denote by $M^{\leq A}$ the closed $A$-neighborhood in
$M$ of $M\ssm\big(\bigcup_{e\in \End(M)} \V_e\big)$.

For every divergent geodesic $\ell$ in $M$, there exist two (possibly
equal) ends $\ell_-$ and $\ell_+$ of $M$ such that
$\lim_{t\ra\pm\infty} \ell(t)=\ell_\pm$. The pair $(\ell_-,\ell_+)$,
which varies in the finite set $\End(M)^2$, is called the {\it type}
of $\ell$. Let $t_-=t_-(\ell)$ be the first time at which $\ell$ exits
the interior of $\V_{\ell_-}$, and $t_+=t_+(\ell)$ be the last time at
which $\ell$ enters the interior of $\V_{\ell_+}$. We define the {\it
  complexity} of the divergent geodesic $\ell$ as
\[
\tau(\ell)= t_+-t_-\geq 0\,.
\]
As usual in counting problems with symmetry, the {\it multiplicity} of
$\ell$ is the inverse of the order of its stabilizer in the orbifold
$M$. See Section \ref{sec:divgeod} for details on multiplicities, and
note that the multiplicities are $1$ when $M$ is a manifold.

We denote by $\Div(M)$ the quotient by the action of $\RR$ by
translation at the source of the set with multiplicities of divergent
geodesics with positive complexity. The set of divergent geodesics
with complexity $0$ is finite up to the action of $\RR$, and can be
ignored in our discussion on equidistribution and counting of
divergent geodesics.  As we want to give counting and equidistribution
results of divergent geodesics with prescribed type, for every nonempty
subset $\T$ of $\End(M)^2$, we consider the space
\begin{equation}\label{eq:defiDivT}
\Div_{\T}(M)=\{\ell\in \Div(M):(\ell_-,\ell_+)\in \T\}
\end{equation}
of divergent geodesics with type contained in $\T$.

We define the {\it Lebesgue measure} $\Leb_\ell$ of $\ell$ in $T^1M$
as the pushforward of the Lebesgue measure of $\RR$ by the map
$t\mapsto \dot\ell(t)=\flow{t}(\dot\ell(0))$.  The study of
equidistribution properties of divergent geodesics is made more
complicated by the fact that $\Leb_\ell$, which is a locally finite
measure with support the geodesic flow orbit $\flow{\RR}(
\dot\ell(0))$, is not a finite measure.

All measures in this paper are locally finite Borel nonnegative
measures. We denote by $\|\mu\|\in[0,+\infty]$ the total mass of a
measure $\mu$ and by $\weakstar$ the weak-star convergence of measures on
locally compact Hausdorff spaces.

Let $\mBM$ be the {\it Bowen-Margulis measure} on $T^1M$, which
is the Liouville measure when $M$ is locally symmetric. We assume that
$\mBM$ is finite, which is for instance the case if $M$ is
locally symmetric.  Then $\mBM$ is mixing for the geodesic flow
by results of Babillot and Dal'Bo, see for instance \cite[\S
  4.2]{BroParPau19}. Its renormalization to a probability measure
$\frac{\mBM}{\|\mBM\|}$ is the unique measure of maximal
entropy on $M$ by \cite{OtaPei04} (and \cite{DilTho25} that removes
the implicit assumption that the sectional curvature of $M$ has
bounded derivative in \cite{OtaPei04}).  The geodesic flow of $M$ has
exponential decay of correlation for the Sobolev regularity with
respect to $\mBM$ for instance if $M$ has constant sectional
curvature by \cite{LiPan22}, or if $M$ is arithmetic locally
symmetric by \cite{KleMar96}, \cite{KleMar99} and \cite{Clozel03}, see
for instance \cite[\S 9.1]{BroParPau19}.

For every end $e\in\End(M)$, let $\sigma^\pm_{e}$ be the (nonzero,
finite) {\it outer/inner skinning measure} on $T^1M$ with support the
outer/inner unit normal bundle of $\partial \V_e$.\footnote{For the
definitions, generalising \cite[\S 1.2]{OhSha12} in constant
curvature, see \cite{ParPau14ETDS} and Section \ref{sec:negcurv}.}
For every nonempty subset $\T$ of $\End(M)^2$, we define a measure on
$T^1M\times T^1M$ by
\begin{equation}\label{eq:defsigmaT}
\sigma_\T=\sum_{(e_-,e_+)\in  \T} \sigma^+_{e_-}\otimes\sigma^-_{e_+}\,.
\end{equation}

\medskip
\noindent\begin{minipage}{5.3cm}
\begin{center}
\includegraphics[width=4.9cm]{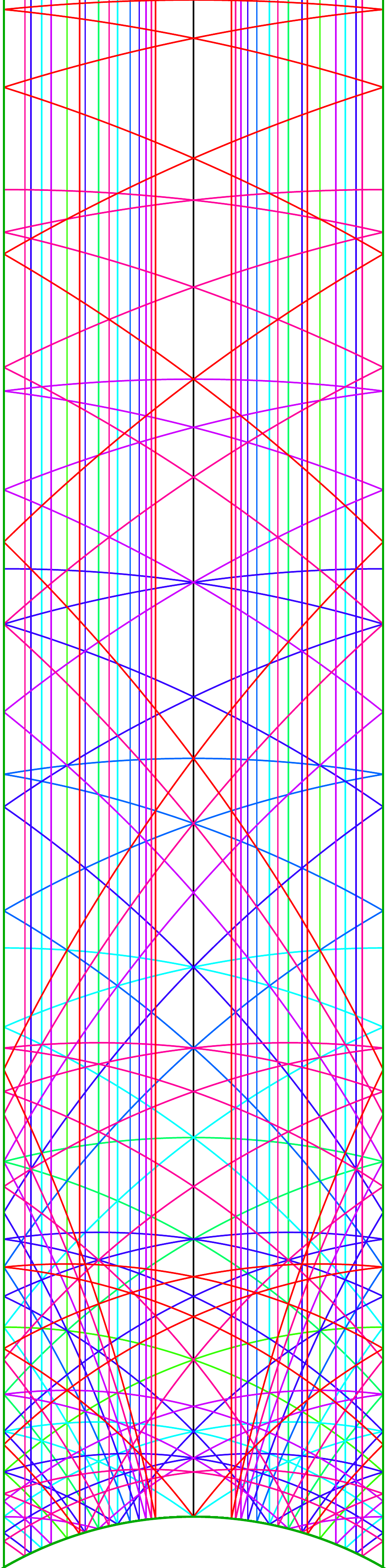}
\end{center}
\end{minipage}
\begin{minipage}{9.6cm}
\setlength\parindent{15pt}
\noindent {\bf Example.}  Let $\hdr$ be the upper halfplane model of
the real hyperbolic plane with constant curvature $-1$. The group
$G=\PSLR$ acts isometrically by homographies on the space $\hdr$, by
the map $(\ga,z)\mapsto \ga\cdot z= \frac{az+b}{cz+d}$ for all
$z\in\hdr$ and $\ga=\pm\begin{pmatrix} a&b\\ c&d\end{pmatrix}\in
G$. Let $\Ga=\PSL_2(\ZZ)$ be the {\it modular group}, which is a
nonuniform arithmetic lattice in $\PSLR$. Let $M=\Ga\bs\hdr$ be the
{\it modular curve}, which is a noncompact complete connected real
hyperbolic good orbifold, with one end. Its standard Margulis cusp
neighbourhood is the $\Ga$-orbit of the horoball $H_\infty$ that
consists of the points $z\in\hdr$ with $\Im \;z\ge 1$.

The figure on the left shows all divergent geodesics of complexity at
most $\ln 10$ for the standard Margulis cusp neighbourhood. These divergent 
geodesics are the images in $M$ of the vertical lines in $\hdr$ with
points at infinity $\infty$ and $\frac pq$, where $p,q\in\ZZ$, $q>0$,
$|p|\le\frac q2$, $1\le q\le 10$. They are represented lifted to the
standard fundamental domain of $\Ga$ 
\[
\Big\{z\in\hdr:-\frac12<\Re\; z<\frac 12,\ |z|>1\Big\}
\]
with side identifications given by $z\mapsto z+1$ and $z\mapsto
-\frac 1z$.

The figure below shows the analogous set of divergent geodesics with
complexity at most $\ln 30$, in the lower part of the fundamental
domain. This illustrates the equidistribution of divergent geodesics
stated in the following Theorem \ref{theo:mainintro}
(\hyperlink{maintheo1}{1}), see also \cite[Theo.~1.5]{DavSha18} in
this specific case.

\begin{center}
\includegraphics[width=9.6cm]{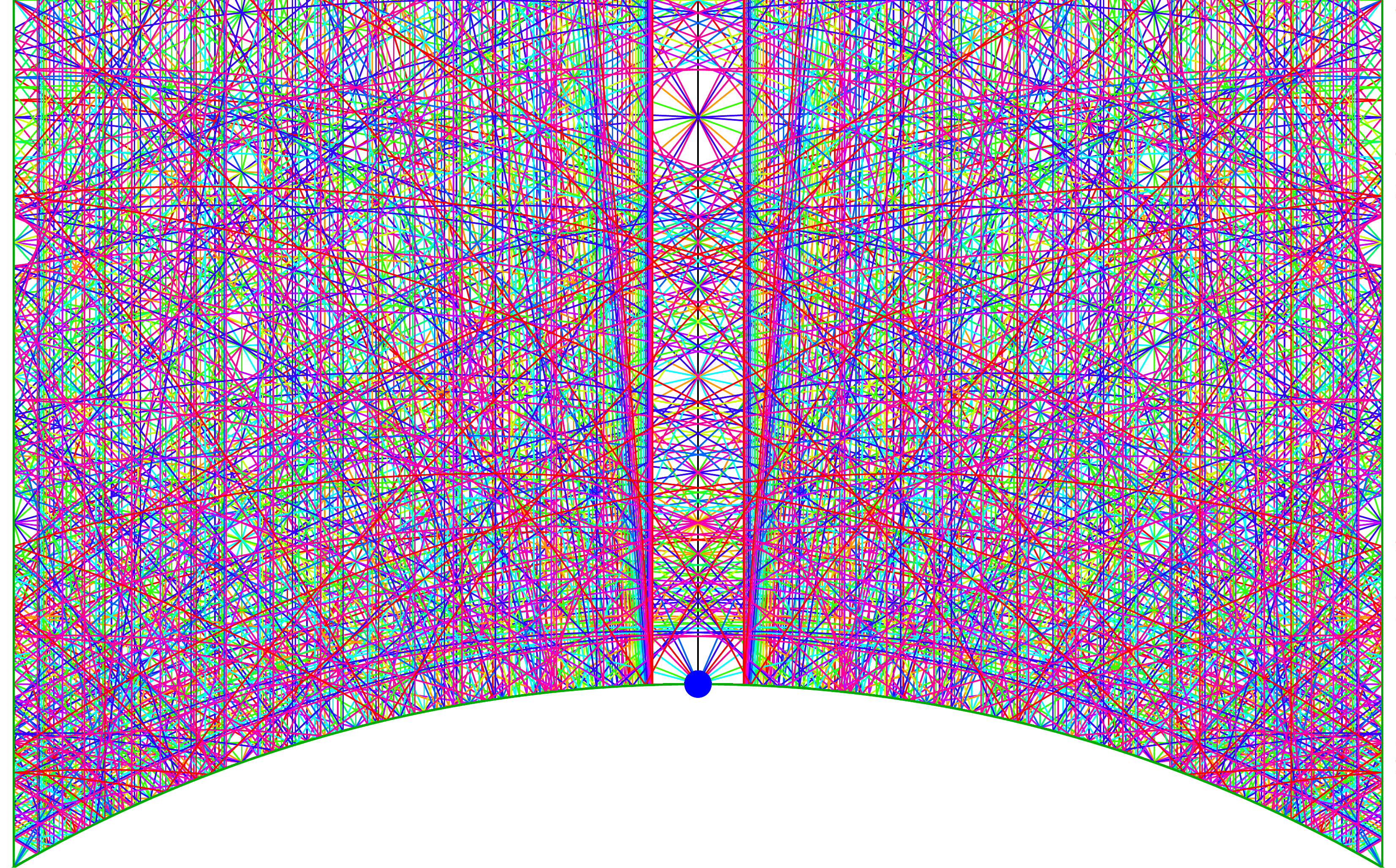}
\end{center}
\end{minipage}

\bigskip
The following is our main result, saying that the divergent geodesics
equidistribute towards the Bowen-Margulis measure as their complexity
tends to $+\infty$.

\btheo \label{theo:mainintro} Let $M$ be a finite volume complete
connected Riemannian good orbifold with dimension at least $2$ and
pinched negative curvature at most $-1$.  Let $\T$ be a nonempty
subset of $\End(M)^2$. Assume that the Bowen-Margulis measure on
$T^1M$ is finite.  

\smallskip\noindent \hypertarget{maintheo1}{(1)} As $T\ra+\infty$, we
have
\[
\frac{h_M\,\|\mBM\|}{\|\sigma_{\T}\|\;T\;e^{h_MT}}
\;\;\sum_{\ell\in\Div_{\T}(M):\;\tau(\ell)\leq T}\Leb_\ell\quad\weakstar\quad
\frac{\mBM}{\|\mBM\|}\;.
\]
If, furthermore, $M$ is locally symmetric with exponential decay of
correlations, then there exists $k\in\NN$ such that for every $\Phi\in
C_c^k(T^1M)$ with support in $\pi^{-1}(M^{\leq A})$, there is an
additive error term in this equidistribution statement when evaluated
on $\Phi$ of the form $\bigO_A\big(\frac{\|\Phi\|_{k,\infty}}{T}
\big)$, where $\|\;\|_{k,\infty}$ is the $W^{k,\infty}$-Sobolev norm.

\smallskip\noindent
\hypertarget{maintheo2}{(2)} As $T\ra+\infty$, we have
\[
\card\big\{\ell\in\Div_{\T}(M):\tau(\ell)\leq T\big\}\;\sim\;
\frac{\|\sigma_\T\|}{h_M\,\|\mBM\|}\;e^{h_MT}\;.
\]
If, furthermore, $M$ is locally symmetric with exponential decay of
correlations, then there exists $\kappa>0$ and an additive error term in
this counting statement of the form $\bigO\big(e^{(h_M-\kappa)T}\big)$.
\etheo

As previously mentioned, one of the main difficulties of this paper is
that the Lebesgue measure associated with a given divergent geodesic
$\ell$ is an infinite measure. We first reduce the study to the
asymptotic distribution of the ``compact cores'' of the divergent
geodesics in Section \ref{sec:divgeod}. We then apply the
equidistribution results proved in \cite{ParPau24a} for manifolds in
Section \ref{sec:proofmanifold}.

In Section \ref{sec:tree}, we give a version of Theorem
\ref{theo:mainintro} for divergent geodesics in geometrically finite
quotients of uniform trees. We refer to Section \ref{sec:negcurv}, the
beginning of Section \ref{sec:tree} and \cite{Lubotzky91,BasLub01,
  Paulin04b,BroParPau19} for the definitions and for background.

\btheo \label{theo:treeintro} Let $X$ be a uniform simplicial tree
without vertices of degree $1$ or $2$ and let $\Ga$ be a finite
covolume geometrically finite discrete subgroup of $\Isom(X)$.  Let
$\T$ be a nonempty subset of $\End(\Ga\bs X)^2$.  Assume that the
length spectrum of $\Ga\bs X$ is $\ZZ$. As $N\ra+\infty$, for the
weak-star convergence of measures on the locally compact space
$\Ga\bs\G X$ of $\Ga$-orbits of geodesic lines in $X$, we have
\[
\frac{(1-e^{-\delta_\Ga})\; \|m_{\rm BM}\|}
{\| \sigma_\T \|\;N\; e^{\delta_\Ga\, N}}
\;\;\sum_{\ell\in\Div_\T(M)\,:\;\tau(\ell)\leq N}\Leb_\ell
\quad\weakstar\quad \frac{m_{\rm BM}}{\|m_{\rm BM}\|}\,.
\]
\etheo

Theorem \ref{theo:treeintro} is part of a more general result, Theorem
\ref{theo:tree} proved in Section \ref{sec:tree}, that covers also the
complementary case where the length spectrum is $2\ZZ$ and gives error
terms for the $\epsilon$-locally constant regularity. We develop an
analog of the equidistribution results of \cite{ParPau24a} for tree
quotients and horoballs in the proof of Theorem \ref{theo:tree}. In
particular, the handling of the error term is much more involved for
trees than in the manifold case.

By \cite{Lubotzky91}, this theorem in particular applies when $X_G$ is
the Bruhat-Tits tree of a rank one simple algebraic group $G$ over a
nonarchimedean local field and $\Ga$ is an arithmetic lattice in
$G$. For instance, let $K=\FF_q((Y^{-1}))$ be the field of formal
Laurent series over $\FF_q$ with indeterminate $Y^{-1}$ and let
$G=\PGL_2(K)$. Then $\Ga=\PGL_2(\FF_q[Y])$ is an arithmetic lattice in
$G$, called the {\it Nagao lattice} (whose length spectrum is
$2\ZZ$). The quotient of the Bruhat-Tits tree $X_G$ by $\Ga$ is then a
geodesic ray, called the {\it modular ray} (see for instance \cite[\S
  15.2]{BroParPau19}) when endowed with its quotient graph of group
structure. In this special case, Theorem \ref{theo:tree} can be
deduced from the case $n=2$ of \cite[Theo.~1.2]{DanPauSay25}, which
gives a stronger equidistribution result in $\Ga\bs G$ (which factors
over $\Ga\bs\G X_G$). Note that $X_G$ is a $(q+1)$-regular tree with
boundary at infinity the projective line $\PP^1(K)=K\cup \{\infty\}$.

\medskip
\noindent {\bf Example.}  When $q=2$, the divergent geodesics in the
modular ray with complexity at most $6$ with respect to the maximal
precisely invariant family of horoballs in $X_G$ are represented in
the following picture (turned horizontally for convenience compared to
the analogous picture in $\PSL_2(\ZZ)\bs \hdr$, all of the divergent
geodesics meant to be pinched vertically to $[0,+\infty[\,$).  For
each shape of image, the pair $(n,m)$ gives the complexity $n$ and the
number $m$ of divergent geodesics in $\Ga\bs X_G$ with this complexity
and shape (they are no longer determined by the shape of their image
in $\Ga\bs X_G$). The divergent geodesics in $\Ga\bs X_G$ are the
images in $\Ga\bs X_G$ of the geodesic lines in $X_G$ starting from
$\infty$ and ending at $\frac{P}{Q}\mod\FF_2[Y]$ with $(P,Q)\in
\FF_2[Y] \times(\FF_2[Y]\ssm \{0\})$. By \cite{Paulin02}, see also
\cite[\S 4.3]{DanPauSay25}, if $\frac{P}{Q}=\cfrac{1}
{a_1+\cfrac{1}{\cdots+ \cfrac{1}{a_k}}}$ is the continued fraction
expansion of $\frac{P}{Q}$ with $a_i\in\FF_2[Y]$ with degree at least
$1$, then the complexity of the associated divergent geodesic is
$n=2\sum_{i=1}^k\deg a_i$, and $m=(q-1)q^{\frac{n}{2}}$ is the number
of choices of the polynomials $a_i$ with a given sequence of degrees
$(\deg a_1,\dots,\deg a_k)$ that defines the shape of the image.

\begin{center}
  \begin{picture}(0,0)%
\includegraphics{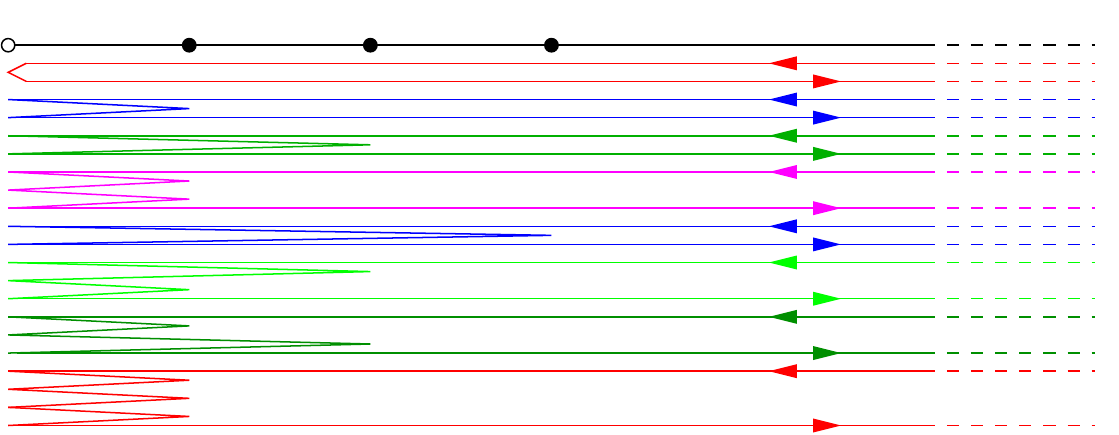}%
\end{picture}%
\setlength{\unitlength}{3812sp}%
\begingroup\makeatletter\ifx\SetFigFont\undefined%
\gdef\SetFigFont#1#2#3#4#5{%
  \reset@font\fontsize{#1}{#2pt}%
  \fontfamily{#3}\fontseries{#4}\fontshape{#5}%
  \selectfont}%
\fi\endgroup%
\begin{picture}(5546,2145)(860,-3523)
\put(901,-1501){\makebox(0,0)[lb]{\smash{{\SetFigFont{11}{13.2}{\rmdefault}{\mddefault}{\updefault}{\color[rgb]{0,0,0}$0$}%
}}}}
\put(1801,-1501){\makebox(0,0)[lb]{\smash{{\SetFigFont{11}{13.2}{\rmdefault}{\mddefault}{\updefault}{\color[rgb]{0,0,0}$1$}%
}}}}
\put(2701,-1501){\makebox(0,0)[lb]{\smash{{\SetFigFont{11}{13.2}{\rmdefault}{\mddefault}{\updefault}{\color[rgb]{0,0,0}$2$}%
}}}}
\put(3601,-1501){\makebox(0,0)[lb]{\smash{{\SetFigFont{11}{13.2}{\rmdefault}{\mddefault}{\updefault}{\color[rgb]{0,0,0}$3$}%
}}}}
\put(6391,-1546){\makebox(0,0)[lb]{\smash{{\SetFigFont{11}{13.2}{\rmdefault}{\mddefault}{\updefault}{\color[rgb]{0,0,0}$\Ga\bs X_G$}%
}}}}
\put(6391,-2581){\makebox(0,0)[lb]{\smash{{\SetFigFont{11}{13.2}{\rmdefault}{\mddefault}{\updefault}{\color[rgb]{0,0,1}$(6,8)$}%
}}}}
\put(6391,-2131){\makebox(0,0)[lb]{\smash{{\SetFigFont{11}{13.2}{\rmdefault}{\mddefault}{\updefault}{\color[rgb]{0,.69,0}$(4,4)$}%
}}}}
\put(6391,-1951){\makebox(0,0)[lb]{\smash{{\SetFigFont{11}{13.2}{\rmdefault}{\mddefault}{\updefault}{\color[rgb]{0,0,1}$(2,2)$}%
}}}}
\put(6391,-1771){\makebox(0,0)[lb]{\smash{{\SetFigFont{11}{13.2}{\rmdefault}{\mddefault}{\updefault}{\color[rgb]{1,0,0}$(0,1)$}%
}}}}
\put(6391,-2806){\makebox(0,0)[lb]{\smash{{\SetFigFont{11}{13.2}{\rmdefault}{\mddefault}{\updefault}{\color[rgb]{0,1,0}$(6,8)$}%
}}}}
\put(6391,-3076){\makebox(0,0)[lb]{\smash{{\SetFigFont{11}{13.2}{\rmdefault}{\mddefault}{\updefault}{\color[rgb]{0,.69,0}$(6,8)$}%
}}}}
\put(6391,-3391){\makebox(0,0)[lb]{\smash{{\SetFigFont{11}{13.2}{\rmdefault}{\mddefault}{\updefault}{\color[rgb]{1,0,0}$(6,8)$}%
}}}}
\put(6391,-2356){\makebox(0,0)[lb]{\smash{{\SetFigFont{11}{13.2}{\rmdefault}{\mddefault}{\updefault}{\color[rgb]{1,0,1}$(4,4)$}%
}}}}
\end{picture}%

\end{center}

\medskip
\noindent{\small {\it Acknowledgements: } The authors thank  the
  French-Finnish CNRS IEA PaCap  for its support.}

\section{Background on negative curvature}
\label{sec:negcurv}

Let $X$ be either a complete simply connected Riemannian manifold $\wt
M$ with dimension at least $2$ and with pinched negative sectional
curvature at most $-1$, or the geometric realisation of a uniform
simplicial tree $\XX$ without vertices of degree $1$ or $2$.  We
denote by $X\cup\partial_\infty X$ the geometric compactification of
$X$. The horoballs in $X$ are closed unless otherwise stated. See for
instance \cite{BriHae99} for background on $\CAT(-1)$ geometry and
discrete groups, and \cite{Serre83,BasLub01}, \cite[\S
  2.6]{BroParPau19} for background on group actions on trees.

Let $x_*\in X$ be a fixed basepoint, with $x_*\in V\XX$ when $X$ is a
tree. For every $\xi\in\partial_{\infty}X$, let $\rho_{\xi}:
[0,+\infty[ \; \ra X$ be the geodesic ray with origin $x_{*}$ and
point at infinity $\xi$. The {\it Busemann cocycle} of $X$ is the
map $\beta: X\times X\times\partial_{\infty} X\to\RR$
defined by
\[
  (x,y,\xi)\mapsto \beta_{\xi}(x,y)= \lim_{t\to+\infty}
  d(\rho_{\xi}(t),x)-d(\rho_{\xi}(t),y)\,.
\]
The {\it visual distance} $d_{x_*}$ on $\partial_\infty X$ seen from
$x_*$ is defined by $d_{x_*}(\xi,\eta)= e^{-\frac{1}{2}(\beta_\xi(x_*,
  \,y) +\beta_\eta(x_*,\,y))}$ where $y$ is the closest point to $x_*$
on the geodesic line $]\xi,\,\xi'[$ between two distinct points at
infinity $\xi$ and $\xi'$.  When $X$ is a tree, we have
\begin{equation}\label{eq:defdistvis}
  d_{x_*}(\xi,\xi')= e^{-d(x_*,\,]\xi,\,\xi'[\,)}\,.
\end{equation}

When $X$ is a manifold, we denote by $\Isom(X)$ the group of
isometries of $\wt M$.  When $X$ is a tree, we denote by $\Isom(X)$
the group of graph automorphisms of $\XX$ without edge inversion.  Let
$\Ga$ be a nonelementary discrete subgroup of $\Isom(X)$. Let $M=\Ga
\bs X$ with its quotient orbifold structure when $X$ is a manifold, and
$M=\Ga \bs X$ with its quotient graph of groups structure when $X$ is
a tree.  When $S$ is a subset or a point in $X$, we denote by $\Ga_S$
the (global) stabiliser of $S$ in $\Ga$.  When $X$ is a tree, we
denote by $V\XX$ the set of vertices of $\XX$, and the {\it covolume}
of $\Ga$ is
\[
\covol \Ga=\sum_{[x]\in \Ga\bs V\XX} \;\frac{1}{\card(\Ga_x)}\,.
\]

Let $\delta_\Ga$ be the critical exponent of $\Ga$. When $X$ is a
manifold and the Bowen-Margulis measure $\mBM$ is finite, $\delta_\Ga$
coincides with the topological entropy $h_M$ of the geodesic flow on
$T^1M$ (see \cite{OtaPei04}):
\begin{equation}\label{eq:criteqentrop}
  \delta_\Ga=h_M\;.
\end{equation}

A continuous mapping $\RR\ra X$ that is an isometric embedding on a
closed interval of $\RR$ and constant, with value in $V\XX$ when $X$
is a tree, on each complementary component is a {\it generalized
  geodesic} in $X$.  We denote by $\gengeod X$ the Bartels-L\"uck
metric space of the generalized geodesic lines $\wt \ell:\RR\ra X$ in
$X$, considering only the generalized geodesic lines $\wt \ell$ such
that $\wt \ell(0)$ is a vertex of $\XX$ when $X$ is a tree. Its
distance is defined, for all $\wt \ell,\wt \ell'\in\gengeod X$ by
\begin{equation}\label{eq:defdistBL}
d(\wt \ell,\wt \ell')=\int_{-\infty}^{+\infty}
d(\wt \ell(t),\wt \ell'(t))\;e^{-2\,|t|}\;dt\;.
\end{equation}
See for instance \cite[\S 2.2]{BroParPau19} for more information.
Geodesic rays in $X$ defined on $\pm\, [0,+\infty[$ are considered as
generalized geodesics by being constant on $\mp\,[0,+\infty[\,$.

We denote the geodesic flow on $\gengeod X$ by $(\flow{t})_{t\in\RR}$
when $X$ is a manifold and by $(\flow{t})_{t\in\ZZ}$ when $X$ is a
tree, with $\flow{t}:\wt\ell\mapsto (s\mapsto \wt\ell(s+t))$. Let $\G
X$ be the subspace of $\gengeod X$ consisting of the geodesic lines of
$X$, which is invariant under the geodesic flow. When $X$ is a
manifold, the space $\G X$ is identified with the unit tangent bundle
$T^1\wt M$ by the $\Isom(X)$-equivariant homeomorphism $\wt\ell\mapsto
\dot{\wt\ell}(0)$.  We denote by $\pi:\gengeod X\ra X$ the footpoint map
$\wt \ell\mapsto \wt \ell(0)$ and again by $\pi: \Ga\bs \gengeod X\ra
\Ga\bs X$ its quotient map. They are proper maps.  For every $\wt
\ell\in\gengeod X$, let $\wt \ell_\pm= \lim_{t\ra\pm\infty} \wt
\ell(t)\in X\cup \partial_\infty X$, which are points at infinity when
$\wt \ell\in\G X$. We denote by $p^-, p^+ : \gengeod X \ra X\cup
\partial_\infty X$ the {\it negative, positive endpoint maps} defined
by $\wt\ell\mapsto \wt\ell_-,\wt\ell_+$ respectively.

For every generalized geodesic $\wt w\in \gengeod X$, which is
isometric on a maximal interval $I$, if $w=\Ga\wt w\in\Ga\bs\gengeod
X$, we define

$\bullet$~ the {\it length}  of $w$ by
\begin{equation}\label{eq:lengthgengeod}
\lambda(w)= \operatorname{length}(I)\in[0,+\infty]\;,
\end{equation}

$\bullet$~ the {\it Lebesgue measure} $\Leb_w$ of $w$ as the measure
on $\Ga\bs\gengeod X$ which is the pushforward by the map $t\mapsto
\flow{t}w$ of the Lebesgue measure on $I$ when $X$ is a manifold, and
of the counting measure on $I\cap \ZZ$ when $X$ is a tree.

For every $w \in \Gamma \bs \gengeod X$, we have $\|\Leb_w\|=
\lambda(w)$.  When $\ell \in\Ga\bs\G X$, we then have $\| \Leb_\ell \|
= +\infty$ and the Lebesgue measure $\Leb_\ell$ is invariant under the
action of the geodesic flow on $\Ga\bs\G X$, since the Lebesgue
measure on $\RR$ and the counting measure on $\ZZ$ are invariant under
translations.

We conclude Section \ref{sec:negcurv} by giving details on the
construction of the Bowen-Margulis measure $\mBM$.  We denote by
$\partial_\infty^2X$ the complement of the diagonal in
$\partial_\infty X \times \partial_\infty X$.  {\it Hopf's
  parametrisation} with respect to the basepoint $x_*$ is the
homeomorphism which identifies $\G X$ with $\partial_\infty^2 X
\times\RR$ when $X$ is a manifold and $\partial_\infty^2 X \times\ZZ$
when $X$ is a tree by the map $\wt\ell\mapsto(\wt\ell_-,\wt\ell_+,s)$,
where $s$ is the signed distance to $\wt\ell(0)$ of the closest point
to $x_*$ on the geodesic line $\wt\ell(\RR)$. Note that a change of
base point only changes the third parameter $s$ by an additive
constant. We fix a Patterson-Sullivan density $(\mu_x)_{x\in X}$ when
$X$ is a manifold (see \cite[\S 4.1]{BroParPau19}), and $(\mu_x)_{x\in
  V\XX}$ when $X$ is a tree (see \cite[\S 4.3]{BroParPau19}), for
$\Ga$ (with zero potential). Since $M$ has finite covolume, the
Patterson-Sullivan measures have full support in $\partial_\infty
X$. The {\it Bowen-Margulis measure} on $\G X$ (associated with this
Patterson-Sullivan density) is the measure $\wtmBM$ on $\G X$ given by
the density
\[
d\,\wtmBM(\wt\ell)=
e^{-\delta_\Ga(\beta_{\wt\ell_-}\,(\wt\ell(0),\,x_*)\,+\,\beta_{\wt\ell_+}(\wt\ell(0),\,x_*))}\;
d\mu_{x_*}(\wt\ell_-)\,d\mu_{x_*}(\wt\ell_+)\,ds
\]
in Hopf's parametrisation of $\G X$ with respect to $x_*$. When $X$ is
a tree, we have
\begin{equation}\label{eq:defiBMtree}
  d\,\wtmBM(\wt\ell)= e^{-2 \,\delta_\Ga\, d(x_*, \; ]\,\wt\ell_-,\, \wt\ell_+[ \;)}
  \,d\mu_{x_*}(\wt\ell_-)\,d\mu_{x_*}(\wt\ell_+)\,ds\,.
\end{equation}
The Bowen-Margulis measure $\wtmBM$ is independent of $x_*$, and it is
invariant under the actions of the group $\Ga$ and of the geodesic
flow. Thus, it defines (see \cite[\S 2.6]{PauPolSha15} for the
branched cover issues) a measure $\mBM$ on $\Ga\bs\G X$ which is invariant
under the quotient geodesic flow, called the {\it Bowen-Margulis
  measure} on $\Ga\bs \G X$.

Let $D$ be a nonempty closed convex subset of $X$, which is the
geometric realisation of a subtree of $\XX$ when $X$ is a tree. We
denote by $\partial^1_\pm D$ the {\it outer/inner normal bundle of
  $D$}, that is the subspace of $\gengeod X$ consisting of the
positive/negative geodesic rays $\rho:\pm[0,+\infty[\ra X$ with
$\rho(0)\in \partial D$, $\rho_\pm\in \partial_\infty X\ssm
\partial_\infty D$ and $\rho(0)$ the closest point to $\rho(\pm t)$ on
$D$ for all $t>0$, see \cite[\S 2.4]{BroParPau19} for details.  We
refer to \cite[\S 3]{ParPau14ETDS} when $X$ is a manifold and to
\cite[Chap.~7]{BroParPau19} in general for more background and for the
basic properties of the following measures. The {\it (outer) skinning
  measure} on $\normalout D$ (associated with the above
Patterson-Sullivan density) is the measure $\wt\sigma^+_D$ on
$\normalout D$ defined, using the positive endpoint homeomorphism
$p_+:\rho\mapsto \rho_+$ from $\normalout D$ to $\partial_\infty X
\ssm\partial_\infty D$, by
\[
d\,\wt\sigma^+_D(\rho) = e^{-\delta_\Ga\,
  \beta_{\rho_+}(\rho(0),\,x_*)}\;d\mu_{x_*}(\rho_{+}) \,.
\]
The {\it (inner) skinning measure} $d\,\wt\sigma^-_D(\rho)=
e^{-\delta_\Ga\,\beta_{\rho_-}(\rho(0),\,x_*)}\;d\mu_{x_*}
(\rho_{-})$ is the similarly defined measure on $\normalin D$.  When
$D=\{x_*\}$, we immediately have
\begin{equation}\label{eq:skinningsingle}
  \forall\;\rho\in \normalpm\{x_*\},\qquad d\,\wt
  \sigma_{\{x_*\}}^\pm(\rho)= d\mu_{\{x_*\}}(\rho_\pm)\,.
\end{equation}
If the family $(\ga D)_{\ga \in\Ga/\Ga_{D}}$ is locally finite in $X$,
we denote by $\sigma^\pm_{\Ga \!D}$ the locally finite measure on
$\Ga\bs\gengeod X$ induced by the $\Ga$-invariant locally finite
measure $\sum_{\ga \in\Ga/\Ga_{D}} \ga_*\,\wt\sigma^\pm_D$ on
$\gengeod X$. The support of $\sigma^\pm_{\Ga \!D}$ is contained in
the image $\normalpm (\Ga \!D)$ of $\normalpm D$ by the map $\gengeod
X\ra \Ga\bs\gengeod X$.

\section{Generalities on divergent geodesics}
\label{sec:divgeod}

We assume from now on that $\Ga$ has finite covolume and is
furthermore geometrically finite when $X$ is a tree. We recall below
all the necessary properties, and we refer to \cite{Paulin04b} for the
definition of geometrical finiteness in the case of trees (implied by
the finite covolume assumption when $X$ is a manifold). Note that by
\cite{BasLub01}, there are many more finite covolume tree lattices
than geometrically finite ones.

The set $\End(M)$ of ends of the locally compact topological space $M$
is finite and discrete, and can be described as follows (see
\cite{Bowditch95} for manifolds and \cite{Paulin04b} for trees).  Let
$\Par_\Ga$ be the countable $\Ga$-invariant set of (bounded) parabolic
fixed points of elements of $\Ga$, that is the set of points $\xi\in
\partial_\infty X$ such that its stabilizer $\Ga_\xi$ acts properly
and cocompactly on $\partial_\infty X\ssm\{\xi\}$. Let us choose for
every $\xi\in\Par_\Ga$ any geodesic ray $t\mapsto \xi_t$ in $X$ with
$\lim_{t\ra+\infty}\xi_t =\xi$. Then we have a bijection (independent
of the previous choices) from $\Ga\bs\Par_\Ga$ to $\End(M)$ which
associates to $\Ga\xi$ the end of $M$ towards which converges
$\Ga\xi_t$ as $t\ra +\infty$. For every end $e\in\End(M)$, we fix a
parabolic fixed point $\wh e\in\Par_\Ga$ such that the above bijection
maps $\Ga\wh e$ to $e$.

There exists a $\Ga$-equivariant family $(H_\xi)_{\xi\in \Par_\Ga}$ of
horoballs $H_\xi$ in $X$ centered at every $\xi \in \Par_\Ga$ such that

$\bullet$~ when $X$ is a tree, for every $\xi\in \Par_\Ga$, the
boundary $\partial H_\xi$ of $H_\xi$ is contained in $V\XX$,

$\bullet$~ the open horoballs $\stackrel{\circ}{H_\xi}=H_\xi\ssm
\partial H_\xi$, which are the interiors of the horoballs $H_\xi$, are
pairwise disjoint as $\xi$ ranges over $\Par_\Ga$,

$\bullet$~ the quotient $\Ga\bs\big(X\ssm\bigcup_{\xi\in \Par_\Ga}
\stackrel{\circ}{H_\xi}\big)$ is compact.

\smallskip\noindent Note that $\Ga_{H_\xi}=\Ga_\xi$ for every $\xi \in
\Par_\Ga$. For every $e\in\End(\Ga\bs X)$, the image $\V_e=\Ga H_{\wh
  e}$ in $M=\Ga\bs X$ of the horoball $H_{\wh e}$ is a neighborhood of
the end $e$, called a {\it Margulis cusp neighborhood} of $e$ if $X$
is a manifold, and a {\it cuspidal ray} with point at infinity $e$ if
$X$ is a tree, with respect to the family $(H_\xi)_{\xi\in \Par_\Ga}$.

According to the definition in Section \ref{sec:negcurv}, the outer
(resp.~inner) unit normal bundle $\normalout\V_e=\Ga\normalout H_{\wh
  e}$ (resp.~$\normalin\V_e=\Ga\normalin H_{\wh e}\,$) of the Margulis
cusp neighborhood $\V_e$ is the subset of elements $\Ga\,\wt \ell
_{\mid\, [0,\infty[}$ (resp.~$\Ga\,\wt\ell_{\mid\, ]-\infty,0]}\,$) in
$\Ga\bs \gengeod X$ where $\wt\ell$ is a geodesic line in $\G X$ with
$\wt\ell_-=\wh{e}$ (resp.~$\wt\ell_+=\wh{e}\,$) and $\wt\ell(0)\in
\partial H_{\wh{e}}$. When $X$ is a manifold and $\partial \V_e$ is a
submanifold, then the map $\Ga\,\wt \ell\mapsto \Ga
\,\dot{\wt\ell}(0)$ identifies $\partial^1_+\V_e$ and
$\partial^1_-\V_e\,$ with the two connected components of the unit
normal bundle of the hypersurface $\partial \V_e$, pointing
respectively outwards and inwards from $\V_e$.

A {\it divergent geodesic} in $\Ga \bs \G X$ is an orbit $\ell=\Ga \wt
\ell \in \Ga \bs \G X$ under the action of $\Ga$ of a geodesic line
$\wt \ell\in\G X$ in $X$ both of whose points at infinity are in
$\Par_\Ga$. This corresponds to the definition in the introduction
whether $X$ is a manifold or not, by the above properties of the
family $(H_\xi)_{\xi\in \Par_\Ga}$, and the fact that a geodesic line
that enters in a horoball either converges to its point at infinity or
goes through its boundary after a finite time. The {\it multiplicity}
of $\ell$ (independent of the choice of $\wt \ell$ and of the action
of the geodesic flow on $\ell$) is
\begin{equation}\label{eq:defmultdivgeod}
m(\ell)=\frac{1}{\card\;\Ga_{\wt \ell(\RR)}}\;.
\end{equation}
We define
\[
\ell_-=\lim_{t=-\infty}\ell(t)\in \End(M)\qquad\text{and}\qquad
\ell_+=\lim_{t=+\infty} \ell(t) \in \End(M),
\]
called respectively the {\it starting end} and {\it terminating end}
of $\ell$. The {\it type} of $\ell$ is the pair of ends $(\ell_-,
\ell_+)$ of $M$. We denote by $\Div(M)$ the set with multiplicities of
orbits under the geodesic flow of the divergent geodesics endowed with
their multiplicities. For every $\ell \in \Div(M)$, the Lebesgue
measure $\Leb_{\ell}$ on $\Ga \bs \G X$ as defined in Section
\ref{sec:negcurv} is locally finite, since $\ell:\RR\ra M$ is a proper
map.

Let $\ell\in\Ga\bs\G X$ be a divergent geodesic, and choose $\wt
\ell\in\G X$ such that $\ell=\Ga\wt \ell$. Let $t_-=t_-(\wt \ell)$ be
the time at which $\wt \ell$ exits $H_{\wt \ell_-}$, and $t_+=t_+(\wt
\ell)$ be the time at which $\wt \ell$ enters $H_{\wt\ell_+}$. The
{\it complexity of $\ell$ with respect to the family} $(H_\xi)_{\xi\in
  \Par_\Ga}$ is
\begin{equation}\label{eq:defitau}
\tau(\ell)=t_+-t_-\;,
\end{equation}
which is nonnegative since the horoballs in $(H_\xi)_{\xi\in
  \Par_\Ga}$ have pairwise disjoint interiors. It is independent of
the choice of $\wt \ell$ and of the action of the geodesic flow on
$\ell$, hence we will again denote by $\tau$ the induced map from
$\Div(M)$ to $[0,+\infty[$. This corresponds to the definition in the
introduction when $X$ is a manifold.

Note that the Lebesgue measure $\Leb_\ell$ of a divergent geodesic is
a locally finite measure on $\Ga\bs\G X$, with support the geodesic
flow orbit of $\ell$, but it is an infinite measure.
    
The {\it compact core} of a divergent geodesic $\ell$ in $M$ is the
locally geodesic segment in $M$ denoted by $\alpha_\ell:[0,\tau(\ell)]
\ra M$ such that $\alpha_\ell(t) = \ell(t_-+t)$ for every $t\in
[0,\tau(\ell)]$.  We identify it with its extension to a generalized
geodesic in $\Ga\bs\gengeod X$ which is constant on $]-\infty,0]$ and
$[\tau(\ell),+\infty[$. Note that for every $t\in\RR$ when $X$ is a
manifold and $t\in\ZZ$ when $X$ is a tree, we have
$\alpha_{\flow{t}\ell} =\alpha_{\ell}$.

For all $A\in[0,+\infty[$ and $\xi\in\Par_\Ga$, let $H_\xi[A]$ be the
horoball contained in the horoball $H_\xi$, whose boundary $\partial
H_\xi[A]$ is at distance $A$ from the boundary $\partial H_\xi$ of
$H_\xi$. Let $\stackrel{\circ}{H_\xi}\![A]$ be the interior of
$H_\xi[A]$. The {\it $A$-thick part of $M$ with respect to the family}
$(H_\xi)_{\xi\in \Par_\Ga}$ is the orbispace
\[
M^{\leq A}=\Ga\bs\Big(X\ssm\bigcup_{\xi\in \Par_\Ga}
\stackrel{\circ}{H_\xi}\![A]\Big)\;.
\]
As $A$ tends to $+\infty$, the injectivity radius at a point $x\in
M\ssm M^{\leq A}$ tends to $0$ when $X$ is a manifold, and the order
of the stabilizer of $x\in M\ssm M^{\leq A}$ tends to $+\infty$ when
$X$ is a tree. Since $\Ga$ is geometrically finite with finite
covolume, the $A$-thick part $M^{\leq A}$ is compact of diameter
$\diam \;M^{\leq A}\leq 2A+ \diam\; M^{\leq 0}$.  For every compact subset
$K$ of $M$, there exists $A\in[0,+\infty[$ such that $K$ is contained
in $M^{\leq A}$. Furthermore, any geodesic ray $\rho:[0,+\infty[\;\ra
X$ in $X$ from a point $\rho(0)$ in $\partial H_\xi$ to the point at
infinity $\lim_{t\ra+\infty} \rho(t) =\xi$ meets $\partial H_\xi[A]$
exactly at the point $\rho(A)$.

The following result relates asymptotically the Lebesgue measure of a
divergent geodesic to the one of its compact core.

\blemm\label{lem:divgeocompcore} For every $A\in[0,+\infty[\,$, for
every $f\in C_c(\Ga\bs\gengeod X)$ with support contained in
$\pi^{-1}(M^{\leq A})$, for every divergent geodesic
$\ell\in\Ga\bs\G X$, we have
\[
\big|\;\Leb_\ell(f)-\Leb_{\alpha_\ell}(f)\;\big|\leq 2\,A\,\|f\|_\infty\;.
\]
\elemm

\dem Let $A,f,\ell$ be as in the statement. Since the above inequality
is invariant under the action of the geodesic flow, we may assume that
$t_-(\ell)=0$.  By the definitions of the geodesic flow and of the
footpoint map, we have $\pi(\flow{t}\ell)= \ell(t)$, for every
$t\in\RR$ when $X$ is a manifold and $t\in\ZZ$ when $X$ is a tree.
Hence when $t\leq 0$, we have $d(\pi(\flow{t}\ell), \ell(0)) =|t|$,
since for every $\xi\in \Par_\Ga$, any geodesic ray from the boundary
of $H_\xi$ to its point at infinity injects isometrically in $M=\Ga\bs
X$.  Therefore, if $t\leq -A$, we have $\pi(\flow{t}\ell) \notin
M^{\leq A}$ and in particular $f(\flow{t}\ell)=0$ by the assumption on
the support of $f$. Similarly, we have $f(\flow{t}\ell) =0$ if $t\geq
t_+(\ell)+A$.

Assume first that $X$ is a manifold. By the definition of the measures
$\Leb_\ell$ and $\Leb_{\alpha_\ell}$, by the definition of the
complexity $\tau(\ell)=t_+(\ell)$, we thus have
\begin{align*}
  \big|\;\Leb_\ell(f)-\Leb_{\alpha_\ell}(f)\;\big|&=
  \Big|\;\int_{-\infty}^{+\infty}f(\flow{t}\ell)\,dt-
  \int_{0}^{\tau(\ell)}f(\flow{t}\ell)\,dt\;\Big|\\&=
  \Big|\;\int_{-A}^{t_+(\ell)+A}f(\flow{t}\ell)\,dt-
  \int_{0}^{t_+(\ell)}f(\flow{t}\ell)\,dt\;\Big|\\&=
  \Big|\;\int_{-A}^{0}f(\flow{t}\ell)\,dt+
  \int_{t_+(\ell)}^{t_+(\ell)+A}f(\flow{t}\ell)\,dt\;\Big|
  \leq 2\,A\,\|f\|_\infty\;.
\end {align*}
Though we won't use it when $X$ is a tree, the same proof works, up to
replacing $A$ by $\lfloor A\rfloor$, which does not change
$\pi^{-1}(M^{\leq A})$, and $\int_{t=a}^b$ by $\sum_{t=a}^b$ for
$a,b\in\ZZ\cup\{\pm\infty\}$.  \cqfd

\section{Divergent geodesics and common perpendiculars}
\label{sec:commperp}

In this section, we give an explicit description of the correspondence
between the divergent geodesics in $M$ and the common perpendiculars
between Margulis cusp neighborhoods or cuspidal rays in $M$.

Let $\T\subset\End(M)^2$ be a nonempty set of types of divergent geodesics,
and let 
\[
\Div^+_\T(M)=\{\ell\in \Div(M):(\ell_-,\ell_+)\in\T\,,\,\tau(\ell)>0\}
\]
Note that the set $\Div^+_\T(M)$ differs from the set $\Div_{\T}(M)$
defined in Equation \eqref{eq:defiDivT} (when $X$ is a manifold, but
the definition is valid when $X$ is a tree) only by a finite subset,
hence has the same asymptotic distribution property.

Let $e_-,e_+\in \End(M)$. Since $\Ga$ is geometrically finite (by the
compactness of $M^{\leq 0}$), the cardinality of the set $F_{e_-,e_+}$
of double classes $\llbracket \ga\rrbracket\in \Ga_{\wh{e_-}} \bs \Ga/
\Ga_{\wh{e_+}}$, such that $H_{\wh{e_-}} \cap \ga H_{\wh{e_+}}$ is
nonempty, is finite.  Let $\llbracket\ga\rrbracket\in \Ga_{\wh{e_-}}
\bs\Ga/ \Ga_{\wh{e_+}} \ssm F_{e_-,e_+}$. In particular, we have
$\wh{e_-} \neq \ga\,\wh{e_+}$.  The {\it multiplicity} of
$\llbracket\ga\rrbracket$ is defined (independently of the choice of
the representative $\ga$ of $\llbracket\ga\rrbracket$) by
\begin{equation}\label{eq:defimultiplicity}
m_{\llbracket\ga\rrbracket}
=\frac{1}{\card (\Ga_{\wh{e_-}} \cap \ga\Ga_{\wh{e_+}}\ga^{-1})}\;.
\end{equation}
Let $\wt\ell_{\llbracket\ga\rrbracket}$ be the unique geodesic line in
$\G X$ such that $(\wt\ell_{\llbracket\ga\rrbracket})_-=\wh{e_-}$,
$(\wt\ell_{\llbracket\ga\rrbracket})_+=\ga\,\wh{e_+}$ and
$\wt\ell_{\llbracket \ga\rrbracket}(0)\in\partial H_{\wh{e_-}}$. This
last condition uniquely defines $\wt\ell_{\llbracket\ga\rrbracket}$ in
its orbit under the geodesic flow. Then
$\ell_{\llbracket\ga\rrbracket} =\Ga\wt\ell_{\llbracket\ga\rrbracket}$
is a divergent geodesic in $M$, that does not depend on the choice of
the representative $\ga$ of $\llbracket\ga\rrbracket$.

We refer to \cite{ParPau17ETDS} and \cite{BroParPau19} for background
on common perpendiculars between properly immersed locally convex
closed subsets of $M$, especially to \cite{BroParPau19} when $X$ is a
tree. The compact core $\alpha_{\ell_{\llbracket\ga\rrbracket}}$ of
$\ell_{\llbracket\ga\rrbracket}$ is a common perpendicular between
$\V_{e_-}$ and $\V_{e_+}$, with positive length $\lambda
(\alpha_{\ell_{\llbracket\ga\rrbracket}})$. Furthermore, any common
perpendicular between two Margulis cusp neighborhoods when $X$ is a
manifold, or two cuspidal rays when $X$ is a tree, arises this way:
Indeed, there is a unique way to extend a common perpendicular between
two disjoint horoballs in $X$ to a geodesic line in $X$ whose
points at infinity are the points at infinity of the two horoballs.
Let
\[
\Perp_\T(M)=\bigsqcup_{(e_-,e_+)\in\T} \big(\Ga_{\wh{e_-}}
\bs\Ga/ \Ga_{\wh{e_+}} \ssm F_{e_-,e_+}\big)\;.
\]
\blemm\label{lem:bijperpdiv} Assume that $X$ is a manifold. The
mapping from $\Perp_\T(M)$ to $\Div^+_\T(M)$ induced by $\llbracket
\ga \rrbracket\mapsto \ell_{\llbracket\ga\rrbracket}$ is a bijection
between sets with multiplicities. For every $\llbracket\ga\rrbracket
\in\Perp_\T(M)$, we have
\begin{equation}\label{eq:samemultlong}
m(\ell_{\llbracket\ga\rrbracket})=m_{\llbracket\ga\rrbracket}
\qquad\text{and}\qquad \tau(\ell_{\llbracket\ga\rrbracket})=
\lambda(\alpha_{\ell_{\llbracket\ga\rrbracket}})\;.
\end{equation}
\elemm

\dem For every $(e_-,e_+)\in\T$, for every divergent geodesic $\ell$
in $M=\Ga\bs X$ with $\ell_\pm=e_\pm$ and $\tau(\ell)>0$, let $\wt
\ell\in\G X$ be such that $\ell=\Ga\wt\ell$. We may assume that
$\wt\ell\,_-=\wh{e_-}$ up to the action of an element of $\Ga$ on $\wt
\ell$ unique modulo left multiplication by an element of
$\Ga_{\wh{e_-}}$, and that $\wt\ell(0)\in \partial H_{\wh{e_-}}$ up to
the action of a unique element of the geodesic flow on $\wt\ell$. Then
there exists $\ga_\ell \in\Ga$ such that $\wt\ell\,_+= \ga_\ell\,
\wh{e_+}\,$, and $\ga_\ell$ is unique modulo multiplication on the
right by an element of $\Ga_{\wh{e_+}}$ and independent of the action
of the geodesic flow on $\wt\ell$.  Since $\tau(\ell)>0$, the
horoballs $H_{\wh{e_-}}$ and $\ga_\ell H_{\wh{e_+}}$ are disjoint.
Hence the double class $\llbracket\ga_\ell\rrbracket\in \Ga_{\wh{e_-}}
\bs\Ga/ \Ga_{\wh{e_+}}$ is well defined and does not belong to
$F_{e_-,e_+}$.

It is clear by construction that the maps $\llbracket\ga\rrbracket
\mapsto \ell_{\llbracket\ga\rrbracket}$ and $\ell \mapsto \llbracket
\ga_\ell \rrbracket$ induce maps from $\Perp_\T(M)$ to $\Div^+_\T(M)$
and from $\Div^+_\T(M)$ to $\Perp_\T(M)$ that are inverse one of the
other. Since $X$ is a manifold (hence nontrivial geodesic segments
have a unique extension to a geodesic line), the stabilizer in $\Ga$
of the common perpendicular between two disjoint horoballs in $X$ is
equal to the stabilizer in $\Ga$ of the geodesic line between the
points at infinity of the two horoballs. The equalities
\eqref{eq:samemultlong} (the first one being the definition of a
bijection between sets with multiplicities) are then immediate by the
definitions of the two multiplicities \eqref{eq:defmultdivgeod} and
\eqref{eq:defimultiplicity}, of the complexity of a divergent geodesic
\eqref{eq:defitau} and of the length of a generalized geodesic
\eqref{eq:lengthgengeod}.  \cqfd

\medskip
For every $T>0$, let us consider the counting function (with
multiplicities) of the common perpendiculars with length at most $T$,
between two Margulis cusp neighborhoods when $X$ is a manifold, or two
cuspidal rays when $X$ is a tree, whose pair of ends belong to $\T$,
defined by (using the definition of sums over sets with
multiplicities for the second equality below)
\[
\N_\T(T)=\sum_{\substack{(e_-,e_+)\in\T\\
\llbracket\ga\rrbracket \in \Ga_{\wh e_-}\bs \Ga/ \Ga_{\wh e_+}\ssm F_{e_-,e_+}\\
\lambda(\alpha_{\ell_{\llbracket\ga\rrbracket}})\leq T}} m_{\llbracket\ga\rrbracket}
=\card\big\{\llbracket\ga\rrbracket \in \Perp_\T(M):
\lambda(\alpha_{\ell_{\llbracket\ga\rrbracket}})\leq T\big\}\;.
\]

For every $(e_-,e_+)\in \T$, we denote respectively by $\sigma^+_{e_-}
= \sigma^+_{\V_{e_-}}$ and $\sigma^-_{e_+}= \sigma^-_{\V_{e_+}}$ the
outer skinning measure of $\V_{e_-}=\Ga H_{\wh{e_-}}$ and inner
skinning measure of $\V_{e_+}=\Ga H_{\wh{e_+}}$ in $\Ga\bs\gengeod X$,
defined at the end of Section \ref{sec:negcurv}.  They are locally
finite measures on $\Ga\bs\gengeod X$ with support contained in
$\partial^1_+\V_{e_-}$ and $\partial^1_-\V_{e_+}$ respectively.  Since
$\Ga_\xi$ acts cocompactly on $\partial H_\xi$ for every
$\xi\in\Par_\Ga$ (hence $\partial^1_\pm \V_e$ is compact for every
$e\in \End(M)$), and since the limit set of $\Ga$ is equal to the
whole $\partial_\infty X$, the skinning measures $\sigma^\pm_{e^\mp}$
are finite and nonzero, with support exactly
$\partial^1_\pm\V_{e^\mp}$.

Assume till the end of Section \ref{sec:commperp} that $X$ is a
manifold and that the Bowen-Margulis measure is finite and mixing
under the geodesic flow, which is the case under the assumptions of
Theorem \ref{theo:mainintro}. By the second claim of
\cite[Coro.~12]{ParPau17ETDS} applied with $\Omega^\pm=(\normalpm
H_{\ga\wh{e_\mp}})_{\ga\in\Ga/\Ga_{\wh{e_\mp}}}$ for every $(e_-,e_+)
\in \T$, and by a finite summation on $(e_-,e_+)\in\T$, we have as
$T\ra+\infty$
\begin{equation}\label{eq:countgrowth}
  \N_\T(T)\sim \sum_{(e_-,e_+)\in\T}
  \frac{\|\sigma^+_{e_-}\|\;\|\sigma^-_{e_+}\|}
  {\delta_\Ga\;\|\mBM\|}\;e^{\delta_\Ga T}\;.
\end{equation}
Furthermore, if $M$ is locally symmetric with exponential decay of
correlation, then by the second claim of \cite[Theo.~15
  (2)]{ParPau17ETDS} (whose assumptions are indeed satisfied, since
$\partial H_\xi$ is smooth for every $\xi\in \Par_\Ga$ and $M$ has
finite volume), there exists $\tau>0$ and an additive error term in
this counting statement of the form $\bigO\big(e^{(\delta_\Ga-\tau)T}\big)$.

By Equation \eqref{eq:defsigmaT}, we have
\begin{equation}\label{eq:mestotsigmaT}
  \|\sigma_{\T}\|=
  \sum_{(e_-,e_+)\in \T} \|\sigma^-_{e_-}\|\;\|\sigma^+_{e_+}\|\,.
\end{equation}
Claim (\hyperlink{maintheo2}{2}) of Theorem \ref{theo:mainintro}
follows from Equation \eqref{eq:criteqentrop}, from Lemma
\ref{lem:bijperpdiv}, and from Equation \eqref{eq:countgrowth} with
its error term.

\section{Equidistribution of divergent geodesics}
\label{sec:proofmanifold}

The aim of this section is to prove Claim (\hyperlink{maintheo1}{1})
of Theorem \ref{theo:mainintro}. We hence assume that $(X,\Ga)$ is as
in Section \ref{sec:negcurv} with $X$ a manifold. The main tool is the
following theorem proved in \cite[Theo.~1]{ParPau24a} applied with
$A^\pm=\V_{e_\pm}$, so that the skinning measures $\sigma^\pm_{e_\mp}
=\sigma^\pm_{\V_{e_\mp}}$ are finite and nonzero. Recall from the
introduction that if the Bowen-Margulis measure $\mBM$ is
finite, then it is mixing for the geodesic flow on $T^1M$, as needed in
loc.~cit..

\btheo\label{theo:peric} If $\mBM$ is finite, then for all
$(e_-,e_+)\in \T$, as $T\ra+\infty$, we have
\[
\frac{\delta_\Ga\,\|\mBM\|}{T\;e^{\delta_\Ga T}}
\;\;\sum_{\llbracket \ga\rrbracket\in\Perp_{\{(e_-,e_-)\}}(M)\,:\;
  \lambda(\alpha_{\ell_{\llbracket\ga\rrbracket}})\leq T}
  \Leb_{\alpha_{\ell_{\llbracket\ga\rrbracket}}}\quad\weakstar\quad
  \|\sigma^+_{{e_-}}\|\;\|\sigma^-_{{e_+}}\|\;\frac{\mBM}{\|\mBM\|}\,.
\]
Furthermore, if $M$ is locally symmetric with exponential decay of
correlation, then there exists $k\in\NN$ such that for every compact
subset $K$ in $T^1M$ and every $C^k$-smooth function $\Psi:
T^1M\ra\RR$ with support in $K$, with $\|\;\|_{k,2}$ the
$W^{k,2}$-Sobolev norm, we have
\begin{align*}
\frac{\delta_\Ga\,\|\mBM\|}{T\;e^{\delta_\Ga T}}
\;\;\sum_{\llbracket \ga\rrbracket\in\Perp_{\{(e_-,e_-)\}}(M)\,:\;
  \lambda(\alpha_{\ell_{\llbracket\ga\rrbracket}})\leq T}
\Leb_{\alpha_{\ell_{\llbracket\ga\rrbracket}}}(\psi)&=
\|\sigma^+_{{e_-}}\|\;\|\sigma^-_{{e_+}}\|\;\frac{\mBM(\psi)}{\|\mBM\|}
\\&\quad +\bigO_K\Big(\frac{\|\psi\|_{k,2}}{T}\Big)\,.
\end{align*}
\etheo

Hence by a finite summation on $(e_-,e_+)\in\T$, by Lemma
\ref{lem:bijperpdiv} and by Equation \eqref{eq:mestotsigmaT}, for
every $f\in C_c(T^1M)$, we have
\begin{equation}\label{eq:appliPeric}
\lim_{T\ra+\infty}\frac{\delta_\Ga\,\|\mBM\|}{T\;e^{\delta_\Ga T}}
\;\;\sum_{\ell\in\Div^+_\T:\;\tau(\ell)\leq T}\Leb_{\alpha_\ell}(f)
=\|\sigma_\T\|\;\frac{\mBM(f)}{\|\mBM\|}\,.
\end{equation}
Furthermore, since $\pi^{-1}(M^{\leq A})$ is compact, with $k$ as in
Theorem \ref{theo:peric}, there is an error term
$\bigO_A\big(\frac{\|f\|_{k,2} }{T} \big)$ when $f\in C^k(T^1M)$ has
support in $\pi^{-1}(M^{\leq A})$ and when $M$ is locally symmetric
with exponential decay of correlations.

\medskip
\noindent {\bf Proof of Theorem \ref{theo:mainintro}
  (\hyperlink{maintheo1}{1}).}  Let $f\in C_c(T^1M)$ be a test
function.  Using sums over sets with multiplicities, by Equation
\eqref{eq:criteqentrop}, by Lemma \ref{lem:divgeocompcore}, by
Equations \eqref{eq:countgrowth} and \eqref{eq:mestotsigmaT}, and
finally by Equation \eqref{eq:appliPeric}, we have
\begin{align}
  &\;\frac{h_M\,\|\mBM\|}{T\;e^{\delta_\Ga T}}
  \;\;\sum_{\ell\in\Div_\T:\;\tau(\ell)\leq T}\Leb_\ell(f)
  \nonumber\\
  =  &\;\frac{\delta_\Ga\,\|\mBM\|}{T\;e^{\delta_\Ga T}}
  \;\;\sum_{\ell\in\Div_\T:\;\tau(\ell)\leq T}\big(\Leb_{\alpha_\ell}(f)
  +\bigO(A\|f\|_\infty)\big)\nonumber\\
    =  &\;\frac{\delta_\Ga\,\|\mBM\|}{T\;e^{\delta_\Ga T}}
    \;\Big(\sum_{\ell\in\Div_\T:\;\tau(\ell)\leq T}\Leb_{\alpha_\ell}(f)\Big)
    +\bigO\Big(\frac{A\|f\|_\infty}T\Big)\label{eq:maintheoconv}\\
    &\!\!\!\!\!\!\!  \underset{ T\to+\infty}\longrightarrow
    \|\sigma_{\T}\|\;\frac{\mBM(f)}{\|\mBM\|}\,.\nonumber
\end{align}
This proves the convergence claim in Theorem \ref{theo:mainintro}
(\hyperlink{maintheo1}{1}).

Now assume till the end of this proof that $M$ is locally symmetric
with exponential decay of correlations.  Let $k\in\NN$ be given by the
error term in Equation \eqref{eq:appliPeric}.  Let us fix $f\in
C^k_c(T^1M)$ and $A\geq 0$ such that the support of $f$ is contained
in $\pi^{-1}(M^{\leq A})$.  Then by Equation \eqref{eq:maintheoconv}
and by the error term in Equation \eqref{eq:appliPeric}, we have
\begin{align*}
  &\;\frac{h_M\,\|\mBM\|}{T\;e^{h_M\,T}}
  \;\;\sum_{\ell\in\Div_\T:\;\tau(\ell)\leq T}\Leb_\ell(f)\\
  =  &\;\|\sigma_{\T}\|\;\frac{\mBM(f)}{\|\mBM\|}+
  \bigO_A\Big(\frac{\|f\|_{k,2}}{T}\Big)+\bigO\Big(\frac{A\|f\|_\infty}T\Big)\\
  =  &\;\|\sigma_{\T}\|\;\frac{\mBM(f)}{\|\mBM\|}+
  \bigO_A\Big(\frac{\|f\|_{k,\infty}}{T}\Big)\,.
\end{align*}
This ends the proof of Theorem \ref{theo:mainintro}
(\hyperlink{maintheo1}{1}). \qed

\medskip
It follows from the arithmeticity results of Margulis and the results
of \cite{LiPan22}, \cite{KleMar96}, \cite{KleMar99} and
\cite{Clozel03} discussed in the introduction that the only case where
the exponential decay of correlations is not known is when $X$ is a
complex hyperbolic space $\HH^n_\CC$ for $n\geq 2$ and $M$ is not
arithmetic.

\section{Equidistribution of divergent geodesics in geometrically finite
  tree quotients}
\label{sec:tree}

We assume in this section that $X$ is the geometric realisation of a
uniform simplicial tree $\XX$ without vertices of degree $1$ or $2$
and that $\Ga$ is a finite covolume geometrically finite nonelementary
discrete subgroup of $\Isom(X)$. Recall from Section \ref{sec:negcurv}
the notation $M=\Ga \bs X$ with its quotient graph of groups
structure, and the notation $\mBM$ for the Bowen-Margulis measure on
$\Ga\bs\G X$. Let $\T$ be a nonempty subset of $\End(M)^2$. We refer
for instance to \cite[\S 2.6]{BroParPau19} for background on trees and
their (discrete-time) geodesic flow, and to \cite[\S 4.4]{BroParPau19}
for background on the Bowen-Margulis measures of their discrete groups
of automorphisms.

The {\it length spectrum} ${\bf L}_\Ga$ of $\Ga$ is the subgroup of
$\ZZ$ generated by the translation lengths in $X$ of the elements of
$\Ga$. It is equal to $\ZZ$ or $2\ZZ$ under the above assumptions on
$(\XX,\Ga)$ (see \cite[Lem.~4.18]{BroParPau19}, using the fact that
since $\Ga$ has finite covolume besides being geometrically finite,
the minimal nonempty $\Ga$-invariant subtree of $\XX$ is equal to
$\XX$).  For instance, when $\Ga$ is the Nagao lattice
$\PGL_2(\FF_q[Y])$ acting on the Bruhat-Tits tree of
$(\PGL_2,\FF_q((Y^{-1}))\,)$ (see \cite[\S 10.2]{BasLub01}, \cite[\S
15.2]{BroParPau19}), then ${\bf L}_\Ga=2\ZZ$ by \cite[\S II.1.2,
Cor.]{Serre83} or \cite[page 331]{BroParPau19}. See \cite[\S
II.2.3]{Serre83} for lots of other examples with length spectrum
$2\ZZ$.

We define $\delta'_\Ga=1-e^{-\delta_\Ga}$, that we will use when ${\bf
L}_\Ga =\ZZ$, and $\delta''_\Ga = 1-e^{-2\delta_\Ga}$, that we will
use when ${\bf L}_\Ga=2\ZZ$.

We choose a fixed vertex $x_*$ of $\XX$ belonging to the
horosphere that bounds one of the horoballs of the family
$(H_\xi)_{\xi\in \Par_\Ga}$. Let $V_{\rm even}\XX$ be the set of
vertices of $\XX$ at even distance from $x_*$ and
\[
\G_{\rm even} X=\{\ell\in\G X : \ell(0)\in V_{\rm even}\XX\}
\quad\text{and}\quad \gengeod_{\rm even} X=
\{\ell\in\gengeod X : \ell(0)\in V_{\rm even}\XX\}\,.
\]
Note that $\G_{\rm even} X$ and $\gengeod_{\rm even} X$ are clopen subsets
of $\G X$ and $\gengeod X$, respectively.

When ${\bf L}_\Ga =2\ZZ$, then the following properties hold.

$\bullet$~ The group $\Ga$ preserves $V_{\rm even} \XX$ and $\G_{\rm
even} X$, as well as their complementary subsets in $V\XX$ and $\G X$.

$\bullet$~ We choose, as we may by the previous point, the
$\Ga$-equivariant family $(H_\xi)_{\xi\in \Par_\Ga}$ so that the
distance between two elements of this family is even.

$\bullet$~ We denote by $m_{\rm BM, \,even}$ the restriction of the
Bowen-Margulis measure $\mBM$ of $\Ga$ to $\Ga \bs \G_{\rm even}
X$. Since the time-one map $\flow{1}$ of the geodesic flow exchanges
$\G_{\rm even} X$ and its complementary subset, we have $\|m_{\rm
BM,\, even}\| = \frac{\|\mBM\|}{2}$.

$\bullet$~ We define $\Div^+_{\T,\,{\rm even}}(M)$ as the set with
multiplicities (given by Equation \eqref{eq:defmultdivgeod}) of the
orbits under the even-time geodesic flow $(\flow{2t})_{t\in \ZZ}$ of
the divergent geodesics $\underline{\ell}\in\Ga\bs \G_{\rm even} X$
such that $(\underline{\ell}\,_-,\underline{\ell}\,_+)\in \T$ and
$\tau(\underline{\ell})>0$.

$\bullet$~ For every $\ell \in \Div^+_{\T,\,{\rm even}}(M)$ and
$\underline{\ell}\in\Ga\bs \G_{\rm even} X$ a representative of
$\ell$, we denote by $\Leb_{\ell,\,{\rm even}}$ the pushforward measure
on $\Ga \bs \G X$ by the map $t\mapsto \flow{2t}\underline{\ell}$ of
the counting measure on $\ZZ$. This measure has support contained in
$\Ga\bs \G_{\rm even} X$. It does not depend on the choice of the
representative $\underline{\ell}$ of $\ell$, by the invariance under
translation of the counting measure on $\ZZ$.

\medskip
We fix $\beta$ and $\varepsilon$ in $]0,1]$. Given a metric space
$(Y,d)$, we denote by $C^\beta_c(Y)$ the normed real vector space of
(uniformly locally) $\beta$-Hölder-continuous functions with compact
support on $Y$, with norm $f\mapsto \|f\|_\beta = \|f\|_\infty
+ \|f\|'_\beta$, using the convention
\begin{equation}\label{eqdefinormhold}
\| f\|'_\beta =\sup_{0<d(x,y)\leq 1}\frac{|f(x)-f(y)|}{d(x,y)^\beta}
\end{equation}
of for instance \cite[\S 3.1]{BroParPau19}. We refer for instance to
the beginning of Chapter 9 in \cite{BroParPau19} for the definition of
the exponential decay of correlation of $\mBM$ or $m_{\rm BM, \,even}$
for the $\beta$-Hölder regularity of the geodesic flow
$(\flow{t})_{t\in\ZZ}$ or $(\flow{2t})_{t\in\ZZ}$.

Since $\Ga\bs\G X$ is a totally disconnected metric space, another
regularity turns out to be useful, as for instance
in \cite{AthGhoPra12,KemPauSha17}. A function $f:\Ga\bs\G X\ra\RR$ is
{\it $\varepsilon$-locally constant} if $f$ is constant on every open
ball of radius $\varepsilon$. We denote by $C_c^{\,\varepsilon {\,\rm
lc}, \,\beta} (\Ga\bs\G X)$ the normed real vector space of compactly
supported $\varepsilon$-locally constant fonctions $f$, with norm
$\|f\|_{\varepsilon {\,\rm lc}, \,\beta}=\epsilon^{-\beta}\|
f\|_\infty$.  By \cite[Rem.~3.11]{BroParPau19}, the space
$C_c^{\,\varepsilon {\,\rm lc}, \,\beta}(\Ga\bs\G X)$ is continuously
included in $C^\beta_c(\Ga\bs\G X)$ thanks to the inequality
\begin{equation}\label{eq:relathomdnormlcnorm}
\|\cdot\|_\beta \leq 3 \, \|\cdot\|_{\varepsilon {\,\rm
lc}, \,\beta}\;.
\end{equation}
This allows us to state an error term in our theorem,
and the proof will juggle between the norm for general
$\beta$-Hölder-continuous functions and the norm
$\|\cdot\|_{\varepsilon {\, \rm lc}, \, \beta}$.

In what follows, $N$ varies in $\NN$.

\btheo \label{theo:tree}
Let $X,\Ga,\T$ be as in the beginning of Section \ref{sec:tree}, and
let $\beta,\varepsilon\in \;]0,1]$.

\smallskip\noindent\hypertarget{theo:tree1}{(1)}
Assume that ${\bf L}_\Ga =\ZZ$. As $N\ra+\infty$, for the weak-star
convergence of measures on the locally compact space $\Ga\bs\G X$, we
have
\begin{equation}\label{eq:theo:tree1}
\frac{\delta'_\Ga\; \|m_{\rm BM}\|}{\|\sigma_\T\|\;N\;e^{\delta_\Ga\, N}}
\;\;\sum_{\ell\in\Div_\T^+(M)\,:\;\tau(\ell)\leq N}\Leb_\ell
\quad\weakstar\quad \frac{m_{\rm BM}}{\|m_{\rm BM}\|}\,,
\end{equation}
with an additive error term of the form
\[
\bigO\Big(\frac{(A+1)\,
\|\Phi\|_{\varepsilon {\,\rm lc},\,\beta}}{N}\Big)
\]
in this equidistribution statement, when evaluated on $\Phi \in
C_c^{\,\varepsilon {\,\rm lc,\,\beta}}(\Ga\bs\G X)$, where $A \geq 0$
is such that the support of $\Phi$ is contained in $\pi^{-1}(X^{\leq
A})$.

\smallskip\noindent\hypertarget{theo:tree2}{(2)}
Assume that ${\bf L}_\Ga =2\ZZ$. As $N\ra+\infty$, for the weak-star
convergence of measures on the locally compact space $\Ga\bs\G_{\rm
even} X$, we have
\[
\frac{\delta''_\Ga\; \|m_{\rm BM,\, even}\|}
{\|\sigma_\T\|\;N\;e^{\delta_\Ga\, N}}\;\;
\sum_{\ell\in\Div^+_{\T,\,{\rm even}}(M)\,:\;\tau(\ell)\leq 2N}\Leb_{\ell, \, \rm even}
\quad\weakstar\quad \frac{m_{\rm BM,\, even}}{\|m_{\rm BM,\, even}\|}\,,
\]
with an additive error term  of the form
\[
\bigO \Big(\frac{(A+1)\,
\|\Phi\|_{\varepsilon {\,\rm lc}, \,\beta}}{N}\Big)
\]
in this equidistribution statement, when evaluated on $\Phi \in
C_c^{\,\varepsilon {\,\rm lc,\,\beta}}(\Ga \bs \G_{\rm even} X)$, where
$A \geq 0$ is such that the support of $\Phi$ is contained in
$\pi^{-1}(X^{\leq A})$.
\etheo

In the coming proof of Theorem \ref{theo:tree}, we will introduce
auxiliary points $x_0\in V\XX$, with $x_0\in V_{\rm even}\XX$ when
$L_\Ga= 2\ZZ$. We will then use the following joint equidistribution
result, extracted from \cite{BroParPau19} as we explain after its
statement. Though expressed in the universal cover $X$, it says that
the images in $M=\Ga\bs X$ of the initial and terminal generalized
geodesics associated with the common perpendiculars between
$\V_{e_\pm}$ and $\Ga x_0$ jointly equidistribute towards their
skinning measures. We first introduce the necessary notation. For
every vertex $x\in V\XX$ and every parabolic fixed point
$\xi\in\Par_\Ga$ such that $x\notin H_\xi$, we denote by
$\wt{\alpha}_{\xi,\,x}$ the shortest geodesic segment in $\gengeod X$
between $H_\xi$ and $x$, starting at time $0$ from $H_\xi$ and
arriving at time $d(H_\xi, \, x)$ at $x$, seen as an element of
$\gengeod X$ being stationary at time $\leq 0$ and at time $\geq
d(H_\xi, \, x)$.  Note that $\flow{d(H_{\xi},\,x)}
\,\wt{\alpha}_{\xi,\,x}$ is then the parametrization of the common
perpendicular between $H_\xi$ and $\{x\}$ arriving at time $0$ at $x$.

\btheo \label{theo:tree_bpp} Let $e_0 \in \End(M)$ and $x_0 \in V\XX$.

\noindent\hypertarget{tree_bpp1}{(1)}
Assume that ${\bf L}_\Ga =\ZZ$. As $N\ra+\infty$, for the weak-star
convergence of measures on the locally compact space $\gengeod
X \times \gengeod X$, we have
\[
\delta'_\Ga\;e^{-\delta_\Ga\,N} \|m_{\rm BM}\| \;
\sum_{\ga\in\Ga \,:\, 0<d(H_{\ga\,\wh{e_0}},\, x_0)\leq N}
\Delta_{\wt{\alpha}_{\wh{e_0},\,\ga^{-1} x_0}} \otimes
\Delta_{\flow{d(H_{\ga\,\wh{e_0}},\, x_0)} \wt{\alpha}_{\ga\,\wh{e_0},\, x_0}}\quad\weakstar
\quad \wt{\sigma}_{H_{\wh{e_0}}}^+ \otimes \;\wt{\sigma}_{\{x_0\}}^-\,.
\]
Furthermore, there exists $\kappa >0$ such that this convergence has an
additive error term, when evaluated on the product function
$(w,w')\mapsto (\phi^-(w),\phi^+(w'))$ of any two Hölder-continuous
functions $\phi^{\pm} \in C^\beta_c(\gengeod X)$, of the form
$\bigO\big( e^{-\kappa\,N}\|\phi^-\|_\beta\|\phi^+\|_\beta\big)$.

\medskip
\noindent\hypertarget{tree_bpp2}{(2)}
Assume that ${\bf L}_\Ga =2\ZZ$ and that $x_0\in V_{\rm even}\XX$. As
$N\ra+\infty$, for the weak-star convergence of measures on the
locally compact space $\gengeod X \times \gengeod X$, we have
\[
\delta''_\Ga\;e^{-2\,\delta_\Ga N} \|m_{\rm BM,\, even}\| 
\sum_{\ga\in\Ga \,:\, 0<d(H_{\ga\,\wh{e_0}},\, x_0) \leq 2N}
\Delta_{\wt{\alpha}_{\wh{e_0},\, \ga^{-1} x_0}} \otimes
\Delta_{\flow{d(H_{\ga\,\wh{e_0}},\, x_0)} \wt{\alpha}_{\ga\,\wh{e_0},\, x_0}} \weakstar
\wt{\sigma}_{H_{\wh{e_0}}}^+ \otimes \,\wt{\sigma}_{\{x_0\}}^-\,.
\]
Furthermore, there exists $\kappa >0$ such that this convergence has an
additive error term, when evaluated on the product function
$(w,w')\mapsto (\phi^-(w),\phi^+(w'))$ of any two Hölder-continuous
functions $\phi^{\pm} \in C^\beta_c(\gengeod X)$, of the form
$\bigO\big(e^{-\kappa\,N} \|\phi^-\|_\beta \|\phi^+\|_\beta \big)$.
\etheo

\dem
Let us explain more precisely where to find these results in the
book \cite{BroParPau19}.  Since $\XX$ is a uniform tree and $\Ga$ has
finite covolume, the Bowen-Margulis measure $\mBM$ is finite
by \cite[Prop.~4.16~(3)]{BroParPau19}. Since $\Ga$ has finite
covolume, the minimal nonempty $\Ga$-invariant subtree of $\XX$ is
$\XX$, hence is uniform without vertices of degree $2$. We recall the
following facts in the two cases of the statement.

\medskip
(\hyperlink{tree_bpp1}{1}) When $L_\Ga=\ZZ$, the Bowen-Margulis
measure $\mBM$ is mixing under the geodesic flow $(\flow{t}
)_{t\in\ZZ}$ by \cite[Prop.~4.17]{BroParPau19} with system of
conductances $\wt c=0$. The convergence part of Assertion
(\hyperlink{tree_bpp1}{1}) then follows
from \cite[Theo.~11.9]{BroParPau19} with system of conductances $\wt
c=0$, applied with $I^-=\Ga/\Ga_{H_{\wh{e_0}}}$, $I^+=\Ga/\Ga_{x_0}$,
$\D^-= (\ga^- H_{\wh{e_0}})_{\ga^- \in I^-}$ and $\D^+=
(\ga^+\{x_0\})_{\ga^+ \in\Ga/\Ga_{x_0}}$.  More precisely, it follows
from the discrete-time zero-potential version of Equation (11.1)
in \cite{BroParPau19} whose summation
gives \cite[Theo.~11.9]{BroParPau19}, applied with $i$ and $j$ the
trivial classes in $\Ga/\Ga_{H_{\wh{e_0}}}$ and $\Ga/\Ga_{x_0}$, and
with a change of variable $\ga\mapsto \ga^{-1}$, since
$\alpha^-_{i,\,\ga^{-1} j}=\wt{\alpha}_{\wh{e_0},\,\ga^{-1} x_0}$ and
$\alpha^+_{\ga\, i, \,j}=\flow{d(H_{\ga\,\wh{e_0}},
x_0)} \wt{\alpha}_{\ga\,\wh{e_0},\,x_0}$ for every $\ga\in\Ga$.

The error term part of Assertion (\hyperlink{tree_bpp1}{1})
follows from \cite[Theo.~12.16]{BroParPau19} with system of
conductances $\wt c=0$ applied with $\DD^-$ the subtree whose
geometric realisation is $H_{\wh{e_0}}$ and $\DD^+=\{x_0\}$, for the
following reasons. 

$\bullet$~ Its hypothesis (1) is satisfied, since $\wh{e_0}$ is a
bounded parabolic fixed point, hence $\Ga_{H_{\wh{e_0}}}$ acts
cocompactly on $\partial H_{\wh{e_0}}$, and $\Ga_{x_0}$ acts (clearly
!)  cocompactly on $\{x_0\}$.

$\bullet$~ These compactness properties imply that the skinning
measures $\sigma^+_{\V_{e_0}}$ and $\sigma^-_{\Ga\{x_0\}}$ are finite.

$\bullet$~ Its hypothesis (2) on the exponential decay of correlations
is satisfied by \cite[Coro.~9.6 (1)]{BroParPau19}, since $\Ga$ is
geometrically finite.

\medskip (\hyperlink{tree_bpp2}{2})
When $L_\Ga=2\ZZ$, the restriction $m_{\rm BM, \,even}$ of the
Bowen-Margulis measure $\mBM$ of $\Ga$ to $\Ga \bs \G_{\rm even} X$ is
mixing under the even-time geodesic flow $(\flow{2t} )_{t\in\ZZ}$
by \cite[Prop.~4.17]{BroParPau19} with system of conductances $\wt
c=0$. Since the basepoint $x_*$ belongs to the boundary of one of the
horoballs of the family $(H_\xi)_{\xi\in\Par_\Ga}$, since the distance
between two elements of the family $(H_\xi)_{\xi\in\Par_\Ga}$ is even,
and since the distance of two points of any given horosphere in a
simplicial tree is even, for every $\xi\in \Par_\Ga$, we have
$\partial H_\xi\subset V_{\rm even}\XX$. In particular, since $x_0\in
V_{\rm even}\XX$, for every $\xi\in\Par_\Ga$ such that $x_0\notin
H_\xi$, the common perpendicular between $H_\xi$ and $\{x_0\}$ has
both endpoints in $V_{\rm even}\XX$.  Furthermore, the skinning
measures $\sigma^+_{\V_{e_0}}$ and $\sigma^-_{\Ga\{x_0\}}$ have
support contained in $\Ga\bs\gengeod_{\rm even} X$.

The convergence part of Assertion (\hyperlink{tree_bpp2}{2}) then
follows from \cite[Eq.~(11.28)]{BroParPau19} (rather a preliminary
version of it without summation, and with a change of variables
$\ga\mapsto \ga^{-1}$, as above) with system of conductances $\wt
c=0$, applied with $I^-=\Ga/\Ga_{H_{\wh{e_0}}}$, $I^+=\Ga/\Ga_{x_0}$,
$\D^-= (\ga^- H_{\wh{e_0}})_{\ga^- \in I^-}$ and $\D^+=
(\ga^+\{x_0\})_{\ga^+ \in\Ga/\Ga_{x_0}}$, using the definition of
$\delta_\Ga''$ and the fact that $\|m_{\rm BM,\, even}\|
= \frac{\|m_{\rm BM}\|}2$.

The error term part of Assertion (\hyperlink{tree_bpp2}{2}) follows
from \cite[Remark (ii) page 281]{BroParPau19} with system of
conductances $\wt c=0$ applied with $\DD^-$ the subtree whose
geometric realisation is $H_{\wh{e_0}}$ and $\DD^+=\{x_0\}$.  
\cqfd

\medskip
Symmetrically, for all $x\in V\XX$ and $\xi\in\Par_\Ga$ such that
$x\notin H_\xi$, we define a generalized geodesic
$\wt{\alpha}_{x,\,\xi}: t\mapsto\wt{\alpha}_{\xi,\,x} (d(H_\xi,\,x)
-t)$ with origin $\wt{\alpha}_{x,\,\xi}(0)=x$.  With the assumptions
of Theorem \ref{theo:tree_bpp} (\hyperlink{tree_bpp1}{1}), as
$N\to+\infty$, we similarly have the following weak-star convergence
of measures, with the same error term,
\begin{equation}\label{eq:tree_bpprev1}
\delta'_\Ga\;e^{-\delta_\Ga N} \|m_{\rm BM}\| \;
\sum_{\ga\in\Ga\,:\;0<d(H_{\ga\wh{e_0}}, \,x_0) \leq N}
\Delta_{\wt{\alpha}_{x_0,\,\ga\wh{e_0}}} \otimes
\Delta_{\flow{d(H_{\ga\wh{e_0}},\,x_0)} \wt{\alpha}_{\ga^{-1} x_0}, \,\wh{e_0}}
\quad\weakstar\quad
\wt{\sigma}_{\{x_0\}}^+ \otimes \,\wt{\sigma}_{H_{\wh{e_0}}}^-\,.
\end{equation}
With the assumptions of Theorem \ref{theo:tree_bpp}
(\hyperlink{tree_bpp2}{2}), as $N\to+\infty$, we similarly have, with
the same error term,
\begin{equation}\label{eq:tree_bpprev2}
\delta''_\Ga\;e^{-2\,\delta_\Ga N} \|m_{\rm BM,\, even}\| 
\sum_{\gamma \in \Ga \,:\; 0<d(H_{\ga\wh{e_0}},\,x_0) \leq 2N}
\Delta_{\wt{\alpha}_{x_0,\,\ga\wh{e_0}}} \otimes
\Delta_{\flow{d(H_{\ga\wh{e_0}},\,x_0)} \wt{\alpha}_{\ga^{-1} x_0}, \,\wh{e_0}}
\weakstar\;\;
\wt{\sigma}_{\{x_0\}}^+ \otimes \,\wt{\sigma}_{H_{\wh{e_0}}}^-\,.
\end{equation}

Our proof of Theorem \ref{theo:tree} using Theorem \ref{theo:tree_bpp}
is motivated by the proof of Theorem \ref{theo:peric}
in \cite[Theo.~1]{ParPau24a} for good Riemannian orbifolds. We will
adapt each step of the main convergence claim therein to the present
tree case.  Some reduction steps simplify here compared
to \cite{ParPau24a}, since we are not considering general convex
subsets, but only horoballs. On the other hand, our computation of the
error terms in Theorem \ref{theo:tree} will be much more involved than
the one in \cite[Theo.~1]{ParPau24a}, leading to a stronger result.

\medskip
\noindent{\bf Proof of Theorem \ref{theo:tree}.}
The Bowen-Margulis measure on $\Ga\bs\G X$ is finite by for
instance \cite[Prop.~4.16 (3)]{BroParPau19}. Since $\End(M)$ is a
finite set, we can assume that $\T$ contains only one element
$(e_-,e_+)$, then sum the convergence results and the error terms over
a general subset $\T$ of the finite set $\End(M)^2$.

\smallskip\noindent(\hyperlink{theo:tree1}{1}) Assume that $L_\Ga = \ZZ$.

\medskip
\noindent{\bf\hypertarget{treeStep1}{Step 1.} }
Let us prove that if the convergence claim \eqref{eq:theo:tree1} is
true when evaluated on $\Phi \in C_c(\Ga\bs\G X)$ with support in
$\pi^{-1}(M^{\leq 0})$, with an error term
$\bigO\big(\frac{\|\Phi\|_{\varepsilon {\,\rm lc},\,\beta}}{N}\big)$
when besides $\Phi \in C_c^{\,\varepsilon {\,\rm lc,\,\beta}}(\Ga\bs\G
X)$, then Assertion (\hyperlink{theo:tree1}{1}) of
Theorem \ref{theo:tree} is true.

\medskip
Let $A\in\NN$, and let $\Phi\in C_c(\Ga\bs\G X)$ be such that the
support of $\Phi$ is contained in the preimage $\pi^{-1}(M^{\leq A})$
of the $A$-thick part of $M$ with respect to the family
$(H_\xi)_{\xi\in \Par_\Ga}$. Let us denote by $\tau_A:\Div(M)\ra
[0,+\infty[$ the complexity function of the divergent geodesics in $M$
now with respect to the family $(H_\xi[A])_{\xi\in \Par_\Ga}$.  Note
that $\tau_A=\tau+2A$, and that $\Phi$ has support in the preimage by
$\pi$ of the $0$-thick part of $M$ with respect to the family
$(H_\xi[A])_{\xi\in \Par_\Ga}$.

Let $\V_{e_\pm}[A]=\Ga H_{\wh{e_\pm}}[A]$ be the cuspidal ray in $M$
with point at infinity $e_\pm$ with respect to the family
$(H_\xi[A])_{\xi\in\Par_\Ga}$. Let $\sigma^\pm_{e_\mp,\,A}
=\sigma^\pm_{\V_{e_\mp}[A]}$ be the outer/inner skinning measure of
$\V_{e_\mp}[A]$ in $\Ga\bs\gengeod X$. Let $\sigma_{\T,\,A}
=\sigma^-_{e_-,\,A}\otimes\,\sigma^+_{e_+,\,A}$. By \cite[Prop.~4
(iii)]{ParPau14ETDS}, since $\V_{e_\pm}$ is the closed
$A$-neighborhood of $\V_{e_\pm}[A]$, we have
$\|\sigma^\pm_{e_\mp,\,A}\| =e^{-\delta_\Ga
A}\,\|\sigma^\pm_{e_\mp} \|$.  Hence since $\T=\{(e_-,e_+)\}$ and by
Equation \eqref{eq:defsigmaT}, we have $\|\sigma_{\T,\,A}\| =
e^{-2\,\delta_\Ga A}\,\|\sigma_{\T}\|$.

Therefore, under the convergence assumption of
Step \hyperlink{treeStep1}{1}, we have
\begin{align*}
&\frac{\delta'_\Ga\; \|m_{\rm BM}\| \, e^{-\delta_\Ga\, N}}
{N\;\| \sigma_\T \|}
\;\;\sum_{\ell\in\Div_\T^+(M)\,:\;\tau(\ell)\leq N}\Leb_\ell(\Phi)\\ =\;&
\frac{\|\sigma_{\T,\,A}\|}{e^{-2\,\delta_\Ga A}\,\|\sigma_{\T}\|}
\frac{N+2A}{N}\;\frac{\delta'_\Ga\;\|m_{\rm BM}\|\, e^{-\delta_\Ga (N+2A)}}
{(N+2A)\;\| \sigma_{\T,\,A} \|}
\;\sum_{\ell\in\Div_\T^+(M)\,:\;\tau_A(\ell)\leq N+2A}\Leb_\ell(\Phi)\\ 
&\underset{N\to+\infty}\longrightarrow
\frac{m_{\rm BM}(\Phi)}{\|m_{\rm BM}\|}\,.
\end{align*}
If futhermore $\Phi\in C_c^{\,\rm lc,\,\beta}(\Ga\bs\G X)$, then under
the error term assumption of Step \hyperlink{treeStep1}{1}
\begin{align*}
&\frac{\delta'_\Ga\; \|m_{\rm BM}\| \, e^{-\delta_\Ga\, N}}
{N\;\| \sigma_\T \|}
\;\;\sum_{\ell\in\Div_\T^+(M)\,:\;\tau(\ell)\leq N}\Leb_\ell(\Phi)\\ =\;&
\frac{N+2A}{N}\Big(\frac{m_{\rm BM}(\Phi)}{\|m_{\rm BM}\|}
+\bigO\big(\frac{\|\Phi\|_{\varepsilon {\,\rm lc},\,\beta}}{N}\big)\Big)
= \frac{m_{\rm BM}(\Phi)}{\|m_{\rm BM}\|}
+\bigO\Big(\frac{(A+1)\,\|\Phi\|_{\varepsilon {\,\rm lc},\,\beta}}{N}\Big)\,.
\end{align*}
This concludes the proof of Step \hyperlink{treeStep1}{1}.

\medskip
\noindent{\bf\hypertarget{treeStep2}{Step 2.} }
Let us now prove the convergence part \eqref{eq:theo:tree1} of
Assertion (\hyperlink{theo:tree1}{1}) of Theorem \ref{theo:tree}.

\medskip
Let $N\in\NN$. Let $x_0\in V\XX\ssm\bigcup_{\xi\in\Par_\Ga} \stackrel
{\circ}{H_\xi}$, so that its image by the canonical projection $X\ra M
= \Ga\bs X$ belongs to $M^{\leq 0}$. We will allow $x_0$ to vary
everywhere in $V\XX\ssm\bigcup_{\xi\in\Par_\Ga} \stackrel{\circ}
{H_\xi}$ at the end of the proof of Step \hyperlink{treeStep2}{2}.

Using Hopf's parametrisation with respect to $x_0$ (see
Section \ref{sec:negcurv}), we define a measure on $\G X$ by
\[
\nu_N = \sum_{\ga^\pm\in\Ga/\Ga_\wh{e_\pm}\;:\;
0<d(H_{\ga^-\wh{e_-}},\, H_{\ga^+\wh{e_+}}) \leq N}
\Delta_{\ga^-{\wh{e_-}}} \otimes \Delta_{\ga^+{\wh{e_+}}} \otimes ds\,.
\]
This measure is independent of $x_0$, hence is $\Ga$-invariant, by the
invariance of the counting measure $ds$ on $\ZZ$ under
translations. It is locally finite by the local finiteness of the
family $(H_\xi)_{\xi\in \Par_\Ga}$. Since $\T=\{(e_-,e_+)\}$, by the
definitions \eqref{eq:defitau} of the complexity $\tau$
and \eqref{eq:defmultdivgeod} of the multiplicities, the induced
measure (see for instance \cite[\S 2.6]{PauPolSha15} for ramified
cover issues) of $\nu_N$ on $\Ga\bs\G X$ is equal to
$\sum_{\ell\in\Div_{\T}\, :\;\tau(\ell)\leq N} \;\Leb_\ell$ (with the
sum convention over sets with multiplicities).

Let $\psi \in C_c(\G X)$ be such that

$\bullet$~ $\psi$  is  nonnegative,

$\bullet$~ $\psi$ has separate variables when considered, via Hopf's
parametrisation with respect to $x_0$, as defined on a subset of the
product space $\partial_\infty X \times \partial_\infty X \times \ZZ$:
that is, there exist two nonnegative continuous functions
$\psi^\pm:\partial_\infty X\ra \RR$ and a nonnegative bounded function
$\psi^0 : \ZZ \to \RR$ with finite support such that
\begin{equation}\label{eq:psi_product_func}
\forall\;(\xi_-, \xi_+, s)\in\partial_\infty^2X\times\ZZ,\quad
\psi(\xi_-, \xi_+, s)=\psi^-(\xi_-) \;\psi^+(\xi_+) \;\psi^0(s)\,,
\end{equation}

$\bullet$~ $\psi$ has support in the compact-open subset $\pi^{-1}(\{x_0\})
= \{\ell\in\G X \,:\; \ell(0)=x_0\}$ of $\G X$. In particular, by the
definition of Hopf's parametrisation with respect to $x_0$, the map
$\psi^0$ is a multiple of the characteristic function ${\mathbbm
1}_{\{0\}}$ of the singleton $\{0\}$. Since we will, at the end of the
proof of Step \hyperlink{treeStep2}{2}, take linear combinations of
such functions $\psi$, we actually assume that $\psi^0={\mathbbm
1}_{\{0\}}$.

For every $k\in\NN\ssm\{0\}$, let
\[
\A_k = \big\{(\ga^-, \ga^+) \in \Ga/\Ga_{\wh{e_-}} \times \Ga/\Ga_{\wh{e_+}}
\,:\, d(H_{\ga^-\wh{e_-}}, H_{\ga^+\wh{e_+}})=k\big\}
\]
and
\begin{equation}\label{eq:defiSk}
S_k = \sum_{(\ga^-, \,\ga^+)\in\A_k}\;
\psi^-(\ga^- \wh{e_-}) \;\psi^+(\ga^+ \wh{e_+})\,.
\end{equation}
We have
\begin{equation}\label{eq:ineq_Ak_nu'}
\nu_N(\psi) = \sum_{k=1}^{N} S_{k}\,.
\end{equation}
Let us subdivide the index set of the sum defining $S_k$ using $x_0$
as an intermediate point: For all $i,j \in \NN$, we define
\begin{align*}
\A_i^- &= \{\ga^-\in\Ga/\Ga_{\wh{e_-}} : d(H_{\ga^-\wh{e_-}},\, x_0)=i\}\,,
\\ \A_j^+ &= \{\ga^+\in\Ga/\Ga_{\wh{e_+}} : d(x_0,\, H_{\ga^+\wh{e_+}})=j\}\,.
\end{align*}
Thanks to our condition on the support of $\psi$, if an index
$(\ga^-,\ga^+)$ contributes to the sum defining $\nu_N(\psi)$, then
$x_0$ belongs to the geodesic line $]\ga^-\wh{e_-}
,\, \ga^+\wh{e_+}[\,$ with points at infinity $\ga^-\wh{e_-}$ and
$\ga^+\wh{e_+}$. Since $x_0\notin \;\stackrel{\circ}
{H}_{\ga^-\wh{e_-}}\cup \stackrel{\circ}{H}_{\ga^+\wh{e_+}}$, the
vertex $x_0$ then belongs to the common perpendicular
$[H_{\ga^-\wh{e_-}}, H_{\ga^+\wh{e_+}}]$ between $H_{\ga^-\wh{e_-}}$
and $H_{\ga^+\wh{e_+}}$. For such an index $(\ga^-,\ga^+)$, we then
have for every $k\in\NN\ssm\{0\}$ the equivalence
\begin{equation}
\big(\exists\;i\in\llbracket 0,k\rrbracket,\quad
\ga^-\in\A_i^- \;\text{and}\;\; \ga^+\in\A_{k-i}^+\big) \iff
(\ga^-,\ga^+)\in\A_k\,,\label{eq:dissec_Ak}
\end{equation}
where $\llbracket 0,k\rrbracket=\interval 0k\cap\ZZ$. 
The convergence part of Assertion (\hyperlink{tree_bpp1}{1}) of
Theorem \ref{theo:tree_bpp} applied with $e_0=e_-$ and an integration
over the first factor when projected to $\Ga\bs\gengeod X$ give, as
$N\to \infty$, the following weak-star convergence of measures on
$\gengeod X$
\begin{equation}\label{eq:conv_skinning_x_0prep}
\frac{\delta_\Ga' \;e^{-\delta_\Ga\,N}\;\|m_{\rm BM}\|}{\|\sigma_{e_-}^+\|}
\sum_{\ga^-\in \;\bigcup_{i=1}^N \A_i^-}
\Delta_{\flow{d(H_{\ga^-\wh{e_-}},\, x_0)} \wt{\alpha}_{\ga^-\wh{e_-},\, x_0}}
\quad\weakstar\quad \wt{\sigma}_{\{x_0\}}^-.
\end{equation}
Let $\rho_{\ga^-\wh{e_-},\, x_0}:\;]-\infty,0]\ra X$ be the negative
geodesic ray starting from the point at infinity $\ga^-\wh{e_-}$ and
ending at $x_0$ at time $0$, seen as a generalized geodesic that is
constant on $[0,+\infty[\,$. It belongs to the inner normal bundle
$\normalin\{x_0\}$ of the convex subset $\{x_0\}$. It coincides with
$\flow{d(H_{\ga^-\wh{e_-}},\, x_0)} \wt{\alpha}_{\ga^-\wh{e_-},\,
x_0}$ on $[-d(H_{\ga^-\wh{e_-}},\, x_0),+\infty[$. Hence, by
Equation \eqref{eq:defdistBL}, we have
\begin{align}
&d(\rho_{\ga^-\wh{e_-},\, x_0},\;\flow{d(H_{\ga^-\wh{e_-}},\, x_0)}
\wt{\alpha}_{\ga^-\wh{e_-},\, x_0})=\int_{-\infty}^{-d(H_{\ga^-\wh{e_-}},\, x_0)}
2\big|t+d(H_{\ga^-\wh{e_-}},\, x_0)\big|\,e^{2t}\;dt\nonumber\\
=\;&e^{-2\,d(H_{\ga^-\wh{e_-}},\, x_0)}\int_{0}^{+\infty}2u\,e^{-2u}\;du
=\frac{1}{2}\;e^{-2\,d(H_{\ga^-\wh{e_-}},\, x_0)}\leq 1\,.
\label{eq:approxflowwtalpharho}
\end{align}
In particular, the distance in $\gengeod X$ between
$\rho_{\ga^-\wh{e_-},\, x_0}$ and $\flow{d(H_{\ga^-\wh{e_-}},\,
x_0)}\wt{\alpha}_{\ga^-\wh{e_-},\, x_0}$ tends to $0$ uniformly in
$\ga^-$ as $d(H_{\ga^-\wh{e_-}},\, x_0)$ tends to $+\infty$. Therefore
Equation \eqref{eq:conv_skinning_x_0prep} gives, as
$N\to+ \infty$, the following weak-star convergence of measures on
$\normalin\{x_0\}$
\begin{equation}\label{eq:conv_skinning_x_0}
\frac{\delta_\Ga' \;e^{-\delta_\Ga\,N}\;\|m_{\rm BM}\|}{\|\sigma_{e_-}^+\|}
\sum_{\ga^-\in \;\bigcup_{i=1}^N \A_i^-}
\Delta_{\rho_{\ga^-\wh{e_-},\, x_0}}
\quad\weakstar\quad \wt{\sigma}_{\{x_0\}}^-\,.
\end{equation}
The negative endpoint map $\wt\ell\mapsto \wt\ell_-$ from $\gengeod X$
to $X\cup \partial_\infty X$ restricts to a homeomorphism
$p^-:\normalin\{x_0\}\ra \partial_\infty X$. Note that by
Equation \eqref{eq:skinningsingle} applied with $x_*=x_0$, we have
$(p^-)_*\wt{\sigma}_{\{x_0\}}^-=\mu_{\{x_0\}}$. Taking the pushforward
of Equation \eqref{eq:conv_skinning_x_0} by $p^-$ and evaluating it on
$\psi^-$, we hence obtain
\begin{equation}\label{eq:conv_patersully_x_0_psi-}
\lim_{N \ra \infty}\;
\frac{\delta_\Ga' \;e^{-\delta_\Ga\,N}\;\|m_{\rm BM}\|}{\|\sigma_{e_-}^+\|}
\sum_{\ga^-\in\;\bigcup_{i=1}^N \A_i^-} \psi^-(\ga^-\wh{e_-}) =
\mu_{x_0}(\psi^-)\,.
\end{equation}
For every $i \in \NN$, we define
\begin{equation}\label{eq:defiai}
a_i = \sum_{\ga^-\in\,\A_i^-} \psi^-(\ga^-\wh{e_-})\,.
\end{equation}
Let $\eta>0$. By Equation \eqref{eq:conv_patersully_x_0_psi-} and
since the sets $\A_i^-$ for $i\in\NN$ are pairwise disjoint, there
exists $i_0=i_0(\eta) \in \NN\ssm\{0\}$ such that if $s\in \NN$
satisfies $s \geq i_0$, then
\begin{equation}\label{eq:ineq_sum_ai}
\frac{\|\sigma_{e_-}^+\|(\mu_{x_0}(\psi^-)-\eta)}{\delta_\Ga' \;\|m_{\rm BM}\|}
\;e^{\delta_\Ga s} \leq \sum_{i=1}^s a_i \leq \frac{\|\sigma_{e_-}^+\|
(\mu_{x_0}(\psi^-)+\eta)}{\delta_\Ga' \;\|m_{\rm BM}\|} \;e^{\delta_\Ga s}\,.
\end{equation}
By a similar argument, using Equation \eqref{eq:tree_bpprev1} with
$e_0=e_+$, taking the difference between the index values $j-1$ and
$j$, integrating on the second factor when projected to
$\Ga\bs\gengeod X$, approximating long geodesic segments with origin
$x_0$ by positive geodesic rays starting from $x_0$, using the
definition $\delta_\Ga'=1-e^{-\delta_\Ga}$, and pushing forward by the
positive endpoint map $p^+:\normalout\{ x_0\}\ra \partial_\infty X$,
we have
\begin{equation}\label{eq:conv_patersully_x_0_psi+}
\lim_{j \to +\infty}\;\;\frac{\|m_{\rm BM}\| \;e^{-\delta_\Ga j}}{\|\sigma_{e_+}^-\|}
\sum_{\ga^+\in\,\A_j^+} \psi^+(\ga^+\wh{e_+}) = \mu_{x_0}(\psi^+)\,.
\end{equation}
Let us define the smooth functions $f_\eta^\pm :
[0,+\infty[\; \to \RR$ by
\[
\forall\;s\in[0,+\infty[\,,\quad f_\eta^\pm(s) =
\frac{\|\sigma_{e_+}^-\|\, (\mu_{x_0}(\psi^+) \pm \eta)}{\|m_{\rm BM}\|}\;
e^{\delta_\Ga (k-s)}\,.
\]
Then, as in order to obtain Equation \eqref{eq:ineq_sum_ai}, there
exists $j_0=j_0(\eta)\in\NN\ssm\{0\}$ such that, for every $j\in\NN$
with $j \geq j_0$, we have the inequalities
\begin{equation}\label{eq:ineq_sum_psi+}
f_\eta^-(k-j) \leq \sum_{\ga^+\in\,\A_j^+} \psi^+(\ga^+\wh{e_+}) \leq
f_\eta^+(k-j)\,.
\end{equation}
We set
\[
C_{1,\eta} = \frac{\|\sigma_{e_-}^+\|\,(\mu_{x_0}(\psi^-) +\eta)}
{\delta_\Ga' \;\|m_{\rm BM} \|} \quad\text{and}\quad C_{2,\eta} =
\frac{\|\sigma_{e_+}^-\|\,(\mu_{x_0}(\psi^+)+\eta)}{\|m_{\rm BM}\|}\,.
\]
Fix $N\in\NN$ with $N \geq \max\{i_0, j_0\}$. Let $k\in\llbracket 0,
N\rrbracket$. By decomposing the sum $S_k$ defined in
Equation \eqref{eq:defiSk} using Equation \eqref{eq:dissec_Ak} and by
using the definition \eqref{eq:defiai}, we obtain
\begin{equation}\label{eq:ineq_S_K_dissec}
S_k = \sum_{i=0}^k \;a_i \sum_{\ga^+\in\,\A_{k-i}^+} \psi^+(\ga^+\wh{e_+})\,.
\end{equation}
We can decompose the above sum $\sum_{i=0}^k$ into
$\sum_{i=0}^{k-j_0}$, where we can use the upper bound in
Equation \eqref{eq:ineq_sum_psi+}, and the sum $\sum_{i=k-j_0 +1}^{k}$
which is bounded by $\bigO_{j_0}(1) \times \bigO(e^{\delta_\Ga k})
= \bigO_{j_0}(e^{\delta_\Ga k})$ using
Equations \eqref{eq:conv_patersully_x_0_psi-}
and \eqref{eq:conv_patersully_x_0_psi+}. Using Abel's summation
formula, we obtain
\begin{align}
S_k & \leq
\sum_{i=0}^{k-j_0} a_i\;f_\eta^+(i) + \bigO_{j_0}(e^{\delta_\Ga k})
\nonumber\\ & = \Big(\sum_{i=0}^{k-j_0} a_i\Big) \,f_\eta^+
(k-j_0) -\sum_{i=0}^{k-j_0-1} \Big( \sum_{l=0}^i a_l \Big)
(f_\eta^+(i+1)-f_\eta^+(i)) + \bigO_{j_0}(e^{\delta_\Ga k})
\nonumber\\ & = \bigO(e^{\delta_\Ga k}) \bigO_{j_0}(1) -
\sum_{i=0}^{k-j_0-1}\big(\sum_{l=0}^i a_l\big) (f_\eta^+(i+1)-
f_\eta^+(i)) +\bigO_{j_0}(e^{\delta_\Ga \,k})\,.\label{eq:abelSk}
\end{align}
Note that for every $i\in\llbracket0,k\rrbracket$, we have
\begin{align*}
f_\eta^+(i+1)-f_\eta^+(i) &= e^{\delta_\Ga (k-i)}\;
\frac{\|\sigma_{e_+}^-\|\, (\mu_{x_0}(\psi^+) + \eta)}{\|m_{\rm BM}\|}\,
(e^{-\delta_\Ga}-1)
\\&=-\,\delta_\Ga'\; C_{2,\eta}\; e^{\delta_\Ga(k-i)}
=\bigO(e^{\delta_\Ga\,k})\,.
\end{align*}
Decomposing the sum $\sum_{i=0}^{k-j_0-1}$ in Equation \eqref{eq:abelSk}
into on the one hand $\sum_{i=0}^{\min\{i_0-1,k-j_0-1\}}$, where we
use the above estimate which is uniform in $i\in\llbracket
0,i_0+1\rrbracket$, and on the other hand $\sum_{i=i_0}^{k-j_0-1}$
(which vanishes if $k<i_0+j_0-1$), where we use the upper bound in
Equation \eqref{eq:ineq_sum_ai}, we obtain the inequality
\begin{align*}
S_k & \leq
\sum_{i=i_0}^{k-j_0-1} \Big(\sum_{l=0}^i a_l\Big)\, \delta_\Ga'\; C_{2,\eta}
e^{\delta_\Ga(k-i)} + \bigO_{i_0, j_0}(e^{\delta_\Ga k})
\\ & \leq \delta_\Ga' \; C_{1,\eta} \; C_{2,\eta} \; k \; e^{\delta_\Ga k} +
\bigO_{i_0, j_0}(e^{\delta_\Ga \,k})\,.
\end{align*}
By Equation \eqref{eq:ineq_Ak_nu'} and a computation of the
arithmetico-geometric sum (see for
instance \cite[Lem.~2.1]{ParPau23BSMF} for a generalization that is
well-adapted for the current purpose), for $N\in\NN$ large enough, we
hence have
\begin{align}
\nu_N(\psi) & \leq \delta_\Ga' \; C_{1,\eta} \; C_{2,\eta}\;
\frac{N\,e^{\delta_\Ga\,N}}{1-e^{-\delta_\Ga}} +
\bigO_{i_0, j_0}(e^{\delta_\Ga\,N})\nonumber
\\ & = \frac{\|\sigma_{e_-}^+\| \;\|\sigma_{e_+}^-\|\,
(\mu_{x_0}(\psi^-)+\eta)\,(\mu_{x_0}(\psi^+)+\eta)\, N\;e^{\delta_\Ga\,N}}
{\delta_\Ga'\;\|\mBM\|^2}  +
\bigO_{i_0, j_0}(e^{\delta_\Ga\,N})\,.\label{eq:upperestilnuN}
\end{align}
By our assumption on the support of $\psi$ and by the
definition \eqref{eq:defiBMtree} of the Bowen-Margulis measure
$\wtmBM$ with $x_*=x_0$, we have
\begin{equation*}\label{eq:bow_marg_simplification_bigOeps}
\wtmBM(\psi) = (\mu_{x_0} \otimes \mu_{x_0} \otimes ds) (\psi)\,.
\end{equation*}
Taking in Equation \eqref{eq:upperestilnuN} the upper limit as $N\to
+\infty$ then letting $\eta \to 0$, we hence obtain
\[
\limsup_{N\to+\infty} \;\frac{\delta_\Ga'\; \|\mBM\|}
{N\; e^{\delta_\Ga  \,N}\;\|\sigma_{e_-}^+\|\;\|\sigma_{e_+}^-\|}
\;\nu_N(\psi) \leq \frac{1}{\|\mBM\|}\;
(\mu_{x_0} \otimes \mu_{x_0} \otimes ds)(\psi)=
\frac{\wtmBM(\psi)}{\|\mBM\|}\,.
\]

A similar argument, this time using the lower bounds in
Equation \eqref{eq:ineq_sum_ai} and \eqref{eq:ineq_sum_psi+}, gives
the inequality
\[
\liminf_{N\to +\infty}\; \frac{\delta_\Ga'\;\|\mBM\|}
{N\;e^{\delta_\Ga\,N}\;\|\sigma_{e_-}^+\|\;\|\sigma_{e_+}^-\|}
\;\nu_N(\psi) \geq \frac{\wtmBM(\psi)}{\|\mBM\|}\,.
\]
Since $\T=\{(e_-,e_+)\}$ and by
Equation \eqref{eq:defsigmaT}, we have
$\|\sigma_{e_-}^+\|\; \|\sigma_{e_+}^-\|=\|\sigma_\T\|$. Thus
\[
\lim_{N\to +\infty}\; \frac{\delta_\Ga'\;\|\mBM\|}
{N\;e^{\delta_\Ga\,N}\;\|\sigma_\T\|}
\;\nu_N(\psi) = \frac{\wtmBM(\psi)}{\|\mBM\|}\,.
\]

Let $p:X\ra \Ga\bs X$ be the canonical projection.  A standard
argument of  covering the support with sets $\pi^{-1}(\{x_0\})$ for
finitely many $x_0 \in V\XX\ssm\bigcup_{\xi\in\Par_\Ga} \stackrel
{\circ}{H_\xi}$ and of uniform approximation by linear combinations of
functions with separate variables gives us the weak-star convergence
of measures on $(p\circ\pi)^{-1}(M^{\leq 0})$
\begin{equation} \label{eq:weakstarup}
\frac{\delta_\Ga' \;\|\mBM\|\;e^{-\delta_\Ga\,N}}{N\;\|\sigma_\T\|}
\;\nu_N\quad\weakstar\quad \frac{\wtmBM}{\|\mBM\|}\,.
\end{equation}
Since the measure $\nu_N$ on $\G X$ induces the measure
$\sum_{\ell\in\Div_{\T}\, :\;\tau(\ell)\leq N} \;\Leb_\ell$ on
$\Ga\bs\G X$ and by the weak-star continuity of taking induced measures
(see \cite[\S 2.6]{PauPolSha15}), this concludes the proof of
Step \hyperlink{treeStep2}{2}.

\medskip
\noindent{\bf\hypertarget{treeStep3}{Step 3.} }
Let us now prove the error term part  of
Assertion (\hyperlink{theo:tree1}{1}) of Theorem \ref{theo:tree}.

Let $\beta,\varepsilon\in\interval[open left]01$, and let $x_0\in
V\XX\ssm\bigcup_{\xi\in\Par_\Ga} \stackrel {\circ}{H_\xi}$. We fix
$\psi$ in $C_c^{\,\varepsilon {\,\rm lc},\,\beta}(\G X)$, where $\G X$
is endowed with the Bartels-Lück distance \eqref{eq:defdistBL}.  We
assume again that $\psi$ has support contained in $\pi^{-1}(\{x_0\})$
and can be written as in Equation \eqref{eq:psi_product_func}, again
with $\psi^0 =\mathbbm{1}_{\{x_0\}}$. We now prove that we have an
error term of the form $\bigO(\frac{\|\psi\|_{\varepsilon\,{\rm
lc},\, \beta}} {N})$ in the weak-star
convergence \eqref{eq:weakstarup} when evaluated on $\psi$. By the
error term part of Step \hyperlink{treeStep1}{1}, and by a similar
lifting and approximation process, this will conclude the proof of
Step \hyperlink{treeStep3}{3}.

Compared to the proof of Step \hyperlink{treeStep2}{2}, only minor
changes are needed. We keep the notation $\nu_N$, $\A_k$, $S_k$,
$\A^-_i$, $\A^+_j$ therein. We may assume that $\psi\neq 0$.

First, let us prove that the function $\psi^-:\partial_ \infty
X\ra\RR$ is $\varepsilon$-locally constant for the visual distance
$d_{x_0}$ defined in Equation \eqref{eq:defdistvis} applied with
$x_*=x_0$.  By symmetry, the same can be done for $\psi^+$.

Since $\psi\neq 0$, let $(\xi_-,\xi_+)\in\partial_\infty^2 X$ be such
that $\psi^-(\xi_-)\psi^+(\xi_+)\neq 0$. Let $\eta\in\partial_\infty
X$ and $\eta'\in \partial_\infty X$ be at visual distance for
$d_{x_0}$ less than $\varepsilon$ one of each other. Let us prove that
$\psi^-(\eta) =\psi^-(\eta')$, which gives that $\psi^-$ is
$\varepsilon$-locally constant for $d_{x_0}$. We may assume that
$\eta\neq\eta'$. Since the functions $\psi^\pm$ are continuous and
$\partial_\infty X$ has no isolated point, we may assume that
$\eta\neq\xi_+$ and $\eta'\neq\xi_+$.

Using Hopf's parametrisation with respect to $x_0$, we consider the
geodesic lines $\ell$ and $\ell'$ in $\G X$ with parameters
$(\eta,\xi_+,0)$ and $(\eta', \xi_+,0)$ respectively. We may assume
that one of these two geodesic lines, say $\ell$, passes through $x_0$
(at time $0$): Otherwise, our assumption on the support of $\psi$
yields $\psi^-(\eta)\psi^+(\xi_+)=\psi^-(\eta')\psi^+(\xi_+)=0$, hence
$\psi^-(\eta)=\psi^-(\eta')$ since $\psi^+(\xi_+)\neq 0$.

\medskip\noindent {\bf Claim. }
Let us prove that we also have $\ell'(0)=x_0$.

\medskip
\dem Since $\ell_+=\ell'_+=\xi_+$, there exists $p\in V\XX$ such that the
intersection of the images of $\ell$ and $\ell'$ is the geodesic ray
$[p,\xi_+[$ starting from $p$ with point at infinity
$\xi_+$.
\begin{center}
\begin{picture}(0,0)%
\includegraphics{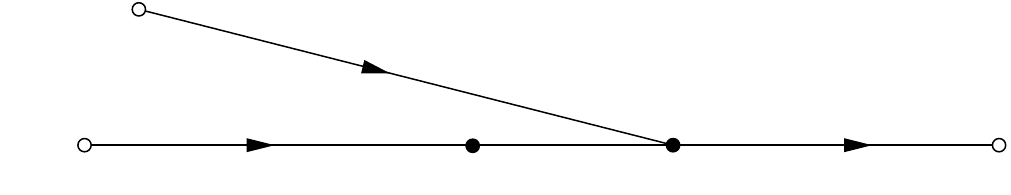}%
\end{picture}%
\setlength{\unitlength}{3812sp}%
\begingroup\makeatletter\ifx\SetFigFont\undefined%
\gdef\SetFigFont#1#2#3#4#5{%
  \reset@font\fontsize{#1}{#2pt}%
  \fontfamily{#3}\fontseries{#4}\fontshape{#5}%
  \selectfont}%
\fi\endgroup%
\begin{picture}(5006,960)(436,-1205)
\put(5221,-1141){\makebox(0,0)[lb]{\smash{{\SetFigFont{11}{13.2}{\rmdefault}{\mddefault}{\updefault}{\color[rgb]{0,0,0}$\ell_+=\ell'_+=\xi_+$}%
}}}}
\put(3781,-1141){\makebox(0,0)[lb]{\smash{{\SetFigFont{11}{13.2}{\rmdefault}{\mddefault}{\updefault}{\color[rgb]{0,0,0}$p$}%
}}}}
\put(2251,-511){\makebox(0,0)[lb]{\smash{{\SetFigFont{11}{13.2}{\rmdefault}{\mddefault}{\updefault}{\color[rgb]{0,0,0}$\ell'$}%
}}}}
\put(1666,-1141){\makebox(0,0)[lb]{\smash{{\SetFigFont{11}{13.2}{\rmdefault}{\mddefault}{\updefault}{\color[rgb]{0,0,0}$\ell$}%
}}}}
\put(2476,-1141){\makebox(0,0)[lb]{\smash{{\SetFigFont{11}{13.2}{\rmdefault}{\mddefault}{\updefault}{\color[rgb]{0,0,0}$\ell(0)=x_0$}%
}}}}
\put(451,-1141){\makebox(0,0)[lb]{\smash{{\SetFigFont{11}{13.2}{\rmdefault}{\mddefault}{\updefault}{\color[rgb]{0,0,0}$\ell_-=\eta$}%
}}}}
\put(451,-421){\makebox(0,0)[lb]{\smash{{\SetFigFont{11}{13.2}{\rmdefault}{\mddefault}{\updefault}{\color[rgb]{0,0,0}$\ell'_-=\eta'$}%
}}}}
\end{picture}%

\end{center}
Assume for a contradiction that $p\in [x_0,\xi_+[\,$, as in the above
picture. Then $x_0$ belongs to the geodesic line $]\eta,\eta'[\,$.
Therefore by Equation \eqref{eq:defdistvis} applied with $x_*=x_0$, we
have $d_{x_0}(\eta,\eta')=e^0=1\geq\varepsilon$. This contradicts
the fact that $d_{x_0}(\eta,\eta')<\varepsilon$.
\begin{center}
\begin{picture}(0,0)%
\includegraphics{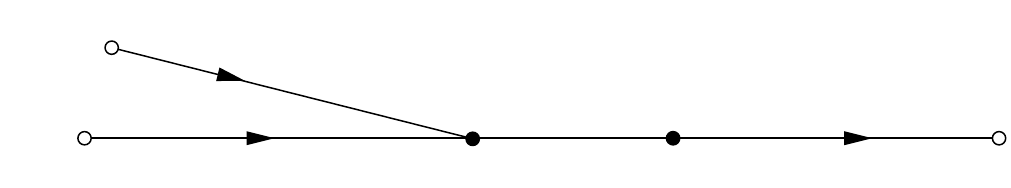}%
\end{picture}%
\setlength{\unitlength}{3812sp}%
\begingroup\makeatletter\ifx\SetFigFont\undefined%
\gdef\SetFigFont#1#2#3#4#5{%
  \reset@font\fontsize{#1}{#2pt}%
  \fontfamily{#3}\fontseries{#4}\fontshape{#5}%
  \selectfont}%
\fi\endgroup%
\begin{picture}(5006,919)(436,-1205)
\put(5221,-1141){\makebox(0,0)[lb]{\smash{{\SetFigFont{11}{13.2}{\rmdefault}{\mddefault}{\updefault}{\color[rgb]{0,0,0}$\ell_+=\ell'_+=\xi_+$}%
}}}}
\put(1666,-1141){\makebox(0,0)[lb]{\smash{{\SetFigFont{11}{13.2}{\rmdefault}{\mddefault}{\updefault}{\color[rgb]{0,0,0}$\ell$}%
}}}}
\put(2746,-1141){\makebox(0,0)[lb]{\smash{{\SetFigFont{11}{13.2}{\rmdefault}{\mddefault}{\updefault}{\color[rgb]{0,0,0}$p$}%
}}}}
\put(3466,-1141){\makebox(0,0)[lb]{\smash{{\SetFigFont{11}{13.2}{\rmdefault}{\mddefault}{\updefault}{\color[rgb]{0,0,0}$\ell(0)=x_0$}%
}}}}
\put(1576,-601){\makebox(0,0)[lb]{\smash{{\SetFigFont{11}{13.2}{\rmdefault}{\mddefault}{\updefault}{\color[rgb]{0,0,0}$\ell'$}%
}}}}
\put(451,-1141){\makebox(0,0)[lb]{\smash{{\SetFigFont{11}{13.2}{\rmdefault}{\mddefault}{\updefault}{\color[rgb]{0,0,0}$\ell_-=\eta$}%
}}}}
\put(451,-421){\makebox(0,0)[lb]{\smash{{\SetFigFont{11}{13.2}{\rmdefault}{\mddefault}{\updefault}{\color[rgb]{0,0,0}$\ell'_-=\eta'$}%
}}}}
\end{picture}%

\end{center}
Hence $p$ belongs to the open geodesic ray $]\eta,x_0[\,$, as in the
above picture. In particular, $x_0$ is the closest point to $x_0$ on
the image of $\ell'$. By the definition of Hopf's parametrisation with
respect to $x_0$, since the third parameter of $\ell'$ is $0$, we have
$\ell'(0)=x_0$, as wanted.
\cqfd

\medskip
With $p$ as in the above proof, let $T\in\NN$ be such that
$p=\ell(-T)=\ell'(-T)$, so that $d(x_0,\;]\eta,\eta'[\,) =T$. By
Equation \eqref{eq:defdistBL}, we have
\[
d(\ell, \, \ell') = \int_{-\infty}^{-T} 2\,|t+T|\;e^{-2\,|t|}\;dt\;=
2\;e^{-2\,T}\int_0^{+\infty} u\;e^{-2\,u}\;du=\frac{1}{2}\;e^{-2\,T}\,.
\]
By Equation \eqref{eq:defdistvis} applied with $x_*=x_0$, we have
$d_{x_0}(\eta,\eta') = e^{-T}$.  Hence
\[
d(\ell, \, \ell') = \frac{1}{2} \;d_{x_0}(\eta,\eta')^2<
\frac{1}{2}\; \varepsilon^2 < \varepsilon\,.
\]
Since $\psi$ is $\varepsilon$-locally constant, this yields
$\psi(\ell)=\psi(\ell')$. Dividing by $\psi_+(\xi_+)\neq 0$, we obtain
$\psi^-(\eta)=\psi^-(\eta')$. Therefore $\psi^-$ is $\varepsilon$-locally
constant, as wanted.

\medskip
As a consequence, since $\psi$ has separate variables, we have
\begin{equation}\label{eq:bound_psi_product}
\|\psi^-\|_{\varepsilon\,{\rm lc},\, \frac{\beta}{2}} \;
\|\psi^+\|_{\varepsilon\,{\rm lc},\, \frac{\beta}{2}}
= \varepsilon^{-\frac{\beta}{2}}\|\psi^-\|_\infty\;
\varepsilon^{-\frac{\beta}{2}}\|\psi^+\|_\infty
= \varepsilon^{-\beta}\|\psi\|_\infty
=\|\psi\|_{\varepsilon\,{\rm lc},\, \beta}\,.
\end{equation}
Let $N\in\NN$. By the error term in Theorem \ref{theo:tree_bpp}
(\hyperlink{tree_bpp1}{1}) for the
$\frac{\beta}{4}$-Hölder-continuity, there exists $\kappa>0$ such that
for every function $\phi^-\in C_c^{\frac{\beta}{4}} (\gengeod X)$,
Equation \eqref{eq:conv_skinning_x_0prep} becomes
\[
\frac{\delta_\Ga' \;e^{-\delta_\Ga \,N}\|\mBM\|}{\|\sigma_{e_-}^+\|}
\sum_{\ga^- \in \,\bigcup_{i=1}^N \A_i^-}
\phi^-(\flow{d(H_{\ga^- \wh{e_-}},\, x_0)}\wt{\alpha}_{\ga^-\wh{e_-},\,x_0})
= \wt{\sigma}_{\{x_0\}}^-(\phi^-) +
\bigO\big(e^{-\kappa N} \|\phi^-\|_{\frac{\beta}{4}}\big)\,.
\]
With $\rho_{\ga^-\wh{e_-},\, x_0}$ as defined after
Equation \eqref{eq:conv_skinning_x_0prep}, since $\phi^-$ is
$\frac{\beta}{4}$-Hölder-continuous and by
Equations \eqref{eqdefinormhold} and \eqref{eq:approxflowwtalpharho},
for all $i\in\NN$ and $\ga^-\in \A_i^-$, we have
\begin{align*}
&\Big|\,\phi^-(\rho_{\ga^-\wh{e_-},\, x_0})-\phi^-(\flow{d(H_{\ga^-\wh{e_-}},\, x_0)}
\wt{\alpha}_{\ga^-\wh{e_-},\, x_0})\,\Big|\\\leq\;\;&
d\big(\rho_{\ga^-\wh{e_-},\, x_0},\;\flow{d(H_{\ga^-\wh{e_-}},\, x_0)}
\wt{\alpha}_{\ga^-\wh{e_-},\, x_0}\big)^{\frac{\beta}{4}}\;\|\phi^-\|'_{\frac{\beta}{4}}
\\\leq\;\;& e^{-\,\frac{\beta}{2}\,d(H_{\ga^-\wh{e_-}},\, x_0)}\;\|\phi^-\|'_{\frac{\beta}{4}}
=e^{-\,\frac{\beta}{2}\,i}\;\|\phi^-\|_{\frac{\beta}{4}}\;.
\end{align*}
Therefore by a geometric series argument and since
$\card\big(\bigcup_{i=1}^N \A_i^-\big)=\bigO(e^{\delta_\Ga \,N})$,
with $\kappa'=\min\{\kappa,\frac{\beta}{2}\}$, for every $\phi^-\in
C^{\frac{\beta}{4}}(\normalin\{x_0\})$,
Equation \eqref{eq:conv_skinning_x_0} becomes
\[
\frac{\delta_\Ga' \;e^{-\delta_\Ga \,N}\|\mBM\|}{\|\sigma_{e_-}^+\|}
\sum_{\ga^- \in \,\bigcup_{i=1}^N \A_i^-}
\phi^-(\rho_{\ga^-\wh{e_-},\, x_0})
= \wt{\sigma}_{\{x_0\}}^-(\phi^-) +
\bigO\big(e^{-\kappa' N} \|\phi^-\|_{\frac{\beta}{4}}\big)\,.
\]
The negative endpoint map $p^- : \normalin\{x_0\} \to \partial_\infty
X$ is $\frac{1}{2}$-Hölder-continuous by a direct adaptation to
negative geodesics rays ending at $x_0$ of the proof of \cite[Lem.~3.4
(4)]{BroParPau19}.  Hence if $\phi^- = \psi^- \circ p^-$, then by
Equation \eqref{eqdefinormhold}, we have $\|\phi^-\|'_{\frac{\beta}
{4}} \leq \|\psi^-\|'_\frac{\beta}{2} \, (\|p^-\|'_\frac{1}{2})
^\frac{\beta} {2}$ and hence $\|\phi^-\|_{\frac{\beta}{4}}
= \bigO(\|\psi^-\|_{\frac{\beta}{2}})$. Applying
the above centered formula to this $\phi^-$,
Equation \eqref{eq:ineq_sum_ai} thus becomes
\begin{equation}\label{eq:estim_sum_ai_errterm}
\forall\;s\in\NN,\qquad
\sum_{i=1}^s \;a_i = \frac{\|\sigma_{e_-}^+\| \; \mu_{x_0}(\psi^-)}
{\delta_\Ga' \;\|\mBM\|} e^{\delta_\Ga \,s} +
\bigO(e^{(\delta_\Ga - \kappa') s} \;\|\psi^-\|_\frac{\beta}{2})\,.
\end{equation}
And similarly, Equation \eqref{eq:ineq_sum_psi+} is replaced by
\begin{equation}\label{eq:estim_sum_psi_+_errterm}
\forall\;j\in\NN,\qquad\sum_{\ga^+\in\A_j^+} \psi^+(\ga^+ \wh{e_+})
= \frac{\|\sigma_{e_+}^-\| \; \mu_{x_0}(\psi^+)}{\|m_{\rm BM}\|}\;
e^{\delta_\Ga \,j} +
\bigO(e^{(\delta_\Ga-\kappa') j} \;\|\psi^+\|_{\frac{\beta}{2}})\,.
\end{equation}
Let us define
\[
C_- = \frac{\|\sigma_{e_-}^+\| \;
\mu_{x_0}(\psi^-)}{\delta_\Ga'\; \|\mBM\|} \quad\text{and}\quad
C_+=\frac{\|\sigma_{e_+}^-\| \; \mu_{x_0}(\psi^+)}{\|\mBM\|}\,.
\]
Using Equation \eqref{eq:relathomdnormlcnorm}, we have
\[
\mu_{x_0}(\psi^\pm) \leq \|\mu_{x_0}\|\;\|\psi^\pm\|_\infty
\leq \|\mu_{x_0}\|\;\| \psi^\pm \|_{\frac{\beta}{2}}
= \bigO(\|\psi^\pm\|_{\varepsilon\,{\rm lc},\, \frac{\beta}{2}})\;.
\]
In particular, we have $C_\pm=\bigO(\|\psi^\pm\|_{\varepsilon\,{\rm
lc},\, \frac{\beta}{2}})$.

Let $k\in\NN\ssm\{0\}$. For every $i\in\llbracket 0, k\rrbracket$, let
$b_i =\sum_{\ga^+\in\A^+_{k-i}} \psi^+ (\ga^+\wh{e_+})$. Note that we
have $b_k=\bigO(\|\psi^+\|_\infty)
= \bigO(\|\psi^+\|_{\varepsilon\,{\rm lc},\, \frac{\beta}{2}})$. By
Equation \eqref{eq:estim_sum_psi_+_errterm}, we have
\begin{equation}\label{eq:estimbi+moinsbi}
b_{i+1}-b_i
=\frac{\|\sigma_{e_+}^-\|\;\mu_{x_0} (\psi^+)}
{\|\mBM\|}\,(e^{-\delta_\Ga}-1)\,e^{\delta_\Ga (k-i)} + \bigO(e^{(\delta_\Ga-\kappa')
(k-i)} \;\|\psi^+\|_{\frac{\beta}{2}})\,.
\end{equation}
Equation \eqref{eq:ineq_S_K_dissec}, Abel's summation formula,
Equations \eqref{eq:estim_sum_ai_errterm}
and \eqref{eq:estimbi+moinsbi}, and the convergence of the series
$\sum_{l \geq 0} e^{-\kappa' \, l}$ yield
\begin{align*}
S_k & = \sum_{i=0}^k \;a_i\,b_i = \Big(\sum_{i=0}^{k}\, a_i\Big) b_k
-\sum_{i=0}^{k-1} \Big( \sum_{l=0}^i \,a_l \Big) (b_{i+1}-b_i) 
\\ & =\big(C_-\,e^{\delta_\Ga \,k}
+\bigO(e^{(\delta_\Ga - \kappa') k} \;\|\psi^-\|_\frac{\beta}{2})\big)
\bigO(\|\psi^+\|_\infty)\\ &\qquad
+\sum_{i=0}^{k-1}\big(C_-\,e^{\delta_\Ga \,i}
+\bigO(e^{(\delta_\Ga - \kappa') i} \;\|\psi^-\|_\frac{\beta}{2})\big)
\big(C_+\,\delta'_k\,e^{\delta_\Ga (k-i)} +\bigO(e^{(\delta_\Ga-\kappa')
(k-i)} \;\|\psi^+\|_{\frac{\beta}{2}})\big)
\\ & = C_- \, C_+ \, \delta_\Ga' \, k \, e^{\delta_\Ga k} +
\bigO(e^{\delta_\Ga k} \|\psi^-\|_{\frac{\beta}{2}}\|\psi^+\|_{\frac{\beta}{2}})\,.
\end{align*}
Hence by Equation \eqref{eq:bound_psi_product}, we have
\[
S_k = C_- \, C_+ \, \delta_\Ga' \, k \, e^{\delta_\Ga k}
+ \bigO(e^{\delta_\Ga k} \|\psi\|_{\varepsilon \, {\rm lc}, \, \beta})\,.
\]
Thus, by the geometric series estimate 
\[
\frac 1{Ne^{\delta_\Ga N}}\sum_{k=0}^Nk\,e^{\delta_\Ga k}
=\sum_{k=0}^N \frac{N-k}{N}
e^{-\delta_\Ga k} = \frac{1}{\delta_\Ga'} + \bigO(\frac{1}{N})\;,
\]
we
have
\[
\frac{\delta_\Ga' \,\|\mBM\|}{\|\sigma_{e_-}^+\|\,
\|\sigma_{e_+}^-\|\, N\, e^{\delta_\Ga \,N}} \, \nu_N(\psi) =
\frac{\wt m_{\rm BM}(\psi)}{\|\mBM\|} +
\bigO\big(\frac{\|\psi\|_{\varepsilon\,{\rm lc},\, \beta}}{N}\big)\,.
\]
This proves that we have an error term of the form
$\bigO(\frac{\|\psi\|_{\varepsilon\,{\rm lc},\, \beta}} {N})$ in the
weak-star convergence \eqref{eq:weakstarup} when evaluated on $\psi$,
hence concluding the proof of Step \hyperlink{treeStep3}{3}, thereby
the one of Assertion (\hyperlink{theo:tree1}{1}) of
Theorem \ref{theo:tree}.

\medskip \noindent(\hyperlink{theo:tree2}{2})
If we assume that $L_\Ga=2\ZZ$, the proof is completely similar to the
one in the case $L_\Ga=\ZZ$ using this time Assertion
(\hyperlink{tree_bpp2}{2}) of Theorem \ref{theo:tree_bpp}.
\cqfd

{\small 
\bibliography{viitteet} 
}

\bigskip
{\small\noindent \begin{tabular}{l} 
J.~P.~and R.~S., Department of Mathematics and Statistics, P.O. Box 35\\ 
40014 University of Jyv\"askyl\"a, FINLAND.\\
{\it jouni.t.parkkonen@jyu.fi, sayousr@jyu.fi}
\end{tabular}

\medskip\noindent \begin{tabular}{l}
F.~P.~and R.~S., Laboratoire de mathématique d'Orsay, UMR 8628 CNRS,\\
Universit\'e Paris-Saclay, 91405 ORSAY Cedex, FRANCE\\
{\it frederic.paulin@universite-paris-saclay.fr,
  rafael.sayous@universite-paris-saclay.fr}
\end{tabular}}

\end{document}